\documentclass{article}
\usepackage{amsmath}
\usepackage{amsfonts,amssymb,cite,mathrsfs,color,authblk}
\date{}

\newtheorem{theorem}{Theorem}[section]
\newtheorem{lemma}{Lemma}[section]

\newtheorem{proposition}{Proposition}[section]

\begin{document}
  \title{ Bounded weak and strong time periodic solutions to a three-dimensional chemotaxis-Stokes model with porous medium diffusion
  \thanks{This work is supported by NSFC(11701384, 11871230), Guangdong Basic and Applied Basic Research Foundation(2020B1515310013).}}

\author[a]{Hailong Ye}
\author[b]{Chunhua Jin\thanks{Corresponding author: jinchhua@126.com }}

\affil[a]{College of Mathematics and Statistics, Shenzhen
University, Shenzhen, 518060, China}
\affil[b]{School of Mathematical Sciences, South China Normal University,
Guangzhou 510631, China}
\renewcommand*{\Affilfont}{\small\it}
\date{}

\maketitle
 \baselineskip=16pt

 \begin{abstract}
In this paper, we study the time periodic problem to a three-dimensional chemotaxis-Stokes model with porous medium diffusion $\Delta n^m$
and inhomogeneous mixed boundary conditions.
By using a double-level approximation method and some iterative techniques, we obtain the existence and time-space uniform boundedness of
weak time periodic solutions for any $m>1$. Moreover, we improve the regularity for $m\le\frac{4}{3}$
and show that the obtained periodic solutions are in fact strong periodic solutions.
\\

\end{abstract}

\bigskip

\section{Introduction and main result}
In this paper, we consider the following coupled
chemotaxis-Stokes model
\begin{equation}\label{problem0}
 \left\{
\begin{array}{l}
n_t+u\cdot\nabla n=\Delta n^m-\chi\nabla\cdot(n\nabla c)+\mu n(a(x,t)-n)+g(x,t), \\[2mm]
c_t+u\cdot\nabla c=\Delta c-cn, \\[2mm]
u_t= \Delta u  -\nabla \pi +n\nabla \varphi, \\[2mm]
\nabla\cdot u =0,
 \end{array} \right.
 \end{equation}
where $m>1$, $(x,t)\in Q=\Omega\times\mathbb{R}^+$, $\Omega\subset\mathbb{R}^3$ is a bounded domain with smooth boundary.
This model describes the motion of oxygen-driven bacteria living in a water drop containing
oxygen. $n, c$ denote the bacterial density, the oxygen concentration respectively,
$n\nabla c$ is the chemotactic flux, $\chi>0$ is the sensitivity coefficient of aggregation induced
by the concentration changes of oxygen, $\mu>0$ is a parameter, $\mu n(a(x,t)-n)$ reflects the proliferation and death
of bacteria in a logistic law, $g(x,t)\ge0$ is source term, $cn$ is the consumption term of oxygen, $u,  \pi$ are the
fluid velocity and the associated pressure, $\nabla \varphi$ is the gravitational
potential. Here, we assume that $a$, $g$ and $\nabla \varphi$ are time periodic functions with
period $T$.

The chemotaxis-fluid model originated from the experimental observation \cite{KE} in concentrated suspensions of swimming bacteria of the species {\it Bacillus subtilis},
and then some biologists and mathematicians began to work on the qualitative description of the pattern formation process in this experiment
 \cite{HP, HPK}. Based on the experimental observations, in 2005, Tuval et al.\cite{Tuval-Cisneros2005}  introduced the following  chemotaxis-fluid model ($m=1$).
\begin{equation}\label{Tuval's problem}
 \left\{
\begin{array}{l}
n_t+u\cdot\nabla n=\Delta n-\nabla\cdot(n\chi(c)\nabla c), \\[1mm]
c_t+u\cdot\nabla c=\Delta c-nf(c), \\[1mm]
u_t+ \tau u\cdot\nabla u= \Delta u  -\nabla \pi +n\nabla \varphi, \\[1mm]
\nabla\cdot u =0,
 \end{array} \right.
 \end{equation}
 which describes the motion of oxygen-driven swimming bacteria in incompressible
fluid, that is, {\it Bacillus subtilis} suspending in a drop of water will swim up an oxygen gradient,
and when the upper bacteria-rich boundary layer is too dense,
it becomes unstable and an overturning instability develops, leading to the formation of falling bacterial
plumes.
In this model, the fluid motion is governed by the Navier-Stokes equations $(\tau=1)$.
However, the viscous force plays a leading role in slow viscous flows,
and the inertial force is far less than the viscous force. Thus, for which,
the Navier-Stokes equations can be approximated using Stokes equations by ignoring
the convective term $u\cdot\nabla u$ (see \cite{Kohr-Pop2004}).

The employed linear diffusion in $\eqref{Tuval's problem}_1$ is the normal type of diffusion associated with
Brownian processes. However, there is evidence that at least some morphogens may
not freely diffuse \cite{Atkins-Paula} and it is needed to develop some nonlinear diffusion models in
biology. Considering that the finite size of bacteria causes the nonlinear enhancement
of random movement of cells at large densities and the diffusion of cells is
more like movement in a porous medium in a viscous fluid, Di Francesco et al. \cite{Francesco-Lorz-Markowich2010}
proposed the following chemotaxis-fluid system by modifying the linear diffusion
$\Delta n$ in \eqref{Tuval's problem} by the nonlinear diffusion $\Delta n^m$:
\begin{equation}\label{nonlinear diffusion problem}
 \left\{
\begin{array}{l}
n_t+u\cdot\nabla n=\Delta n^m-\nabla\cdot(n\chi(c)\nabla c), \\[1mm]
c_t+u\cdot\nabla c=\Delta c-nf(c), \\[1mm]
u_t+ \tau u\cdot\nabla u= \Delta u  -\nabla \pi +n\nabla \varphi, \\[1mm]
\nabla\cdot u =0.
 \end{array} \right.
 \end{equation}
Before going into our mathematical analysis, let us briefly recall some important progresses on system \eqref{nonlinear diffusion problem} and its variants.
For the two dimensional case of \eqref{nonlinear diffusion problem}, that is, $\tau=1$, the global solvability and boundedness
of weak solutions are established completely for any $m>1$ in \cite{Tao-Winkler2012}. While in three dimensional space, the study
of \eqref{nonlinear diffusion problem} with $\tau=0$ is rather tortuous.
In 2011, a global very weak solution $(n\ln n\in L^1)$ is obtained for $m=\frac{4}{3}$ by Liu and Lorz \cite{Liu-Lorz2011}
in dimension 3. Subsequently, Duan and Xiang improved this results to any $m>1$ \cite{Duan-Xiang2014}.
However, this kind of weak solutions may be unbounded, and it is impossible to identify the singularity of the solution.
In 2010, Di Francesco  et al. \cite{Francesco-Lorz-Markowich2010} obtained the existence of a global bounded weak solution for $m\in(m^*, 2]$
($m^*\thickapprox1.81$);
a locally bounded global weak solution was then obtained for $m>\frac{8}{7}$ in 2013 \cite{Tao-Winkler2013}; the uniform boundedness of solutions was
subsequently supplemented for $m>\frac{7}{6}$ \cite{Winkler2015}; further extension was made by
Winkler for $m>\frac{9}{8}$ in \cite{Winkler2018} to a convex domain; recently, Jin \cite{Jin2021arXiv}
obtained the existence of bounded weak solutions for $m>1$ and bounded strong solutions for $m<\frac{5}{4}$.
If the logistic growth term $\mu n(1-n)$ is added to this model,
Lankeit \cite{Lankeit2016} established the global weak solution for the linear diffusion case $m=1$,
and proved that after some waiting time the weak solution becomes smooth and finally converge to
the semi-trivial steady state $(1, 0, 0)$; for the nonlinear diffusion case, Jin \cite{Jin2017JDE} established the
existence of global bounded weak solutions for any $m>1$ to the fluid-free case.

For the time periodic problem of the chemotaxis-fluid models, there are few works concerned.
In 2019, Jin \cite{Jin2020PRSE} considered the chemotaxis-fluid model \eqref{problem0} without source term
for the linear diffusion case ($m=1$),
and proved the existence of bounded strong (and classical) time periodic solutions in dimension 3.
Recently, for the following chemotaxis-Stokes model with porous medium diffusion in dimension 3,
\begin{equation}\label{introduction-problem}
 \left\{
\begin{array}{l}
n_t+u\cdot\nabla n=\Delta n^m-\chi\nabla\cdot(n\nabla c)+\mu n(1-n)+g(x,t), \\[1mm]
c_t+u\cdot\nabla c=\Delta c-c+n, \\[1mm]
u_t= \Delta u  -\nabla \pi +n\nabla \varphi, \\[1mm]
\nabla\cdot u =0,\\[1mm]
\frac{\partial n^m}{\partial \nu}|_{\partial\Omega}=\frac{\partial c}{\partial \nu}|_{\partial\Omega}=u|_{\partial\Omega}=0,
 \end{array} \right.
 \end{equation}
Huang and Jin \cite{Jin2020DCDS} established the existence of uniformly bounded time
periodic solution for any $m\ge\frac{6}{5}$.

It is worth noting that most of the results have been carried out for the closed system,
that is, there is no flux of oxygen and cells through the fluid-air interface.
However, in fact, the experiment is to place the well mixed suspension of {\it Bacillus subtilis} in an open chamber.
On the surface of the water layer, oxygen is allowed to exchange with the outside air. Therefore, for the oxygen in the model,
Dirichlet boundary condition in \cite{CF, LK} or Robin boundary condition in \cite{Braukhoff2017,Wu-Xiang2020} are more realistic.
In \cite{CF, LK}, some numerical results are given, the global existence of small strong solutions around a equilibrium state
is established by Peng and Xiang \cite{PX} in dimension 3, global classical solutions \cite{Braukhoff2017, Braukhoff2020} and
time periodic solutions \cite{Jin2020PRSE} are established respectively in dimension 2 or 3.
Inspired by the above works, we assume that the water drop is surrounded
by air,  oxygen exchange will take place on the boundary of $\Omega$, that is, the solved
oxygen in the water drop may leave, and the free oxygen in the air may diffuse into
the drop. The behaviour of the oxygen exchange can be described by Raoult's law,
which connects the rate of incoming oxygen to the partial vapour pressure of the
oxygen in the surroundings. We assume that the vapour pressure of the free oxygen
is given, and thereby, the incoming rate of oxygen is known. The leaving rate of the
oxygen molecules is proportional to the total number of molecules on the surface.
Therefore, we have the following Robin boundary condition
$$
\frac{\partial c}{\partial \nu}|_{\partial\Omega}=-a_1(x,t)c(x,t)+a_2(x,t),
$$
where $a_1, a_2\in C^{\infty}(\partial\Omega\times[0,\infty))$, $a_1>0$ is the leaving rate
of the oxygen molecules, $a_2\ge0$ with $a_2\not\equiv0$ is the incoming oxygen and depends
on the known vapour pressure of the free oxygen.
By \cite{Braukhoff2017,Lions-Magenes1968}, there exist
\begin{equation}\label{condition1}
g_1, g_2\in C_T^{\infty}(\overline{\Omega}\times[0,\infty))
 \end{equation}
such that
\begin{equation}\label{condition2}
\frac{\partial g_1(x,t)}{\partial \nu}=-a_1(x,t)<0,\quad g_2(x,t)=\frac{a_2}{a_1}\ge0,\quad \frac{\partial g_2(x,t)}{\partial \nu}=0
 \end{equation}
for $(x,t)\in \partial\Omega\times[0,\infty)$. So the boundary condition of $c$ can be rewritened as
$$
\frac{\partial c}{\partial \nu}|_{\partial\Omega}=\frac{\partial g_1(x,t)}{\partial \nu}(c(x,t)-g_2(x,t)).
$$
For $u$, we still consider the no-slip
boundary condition, namely,
\begin{equation}\label{Dirichlet boundary condition for u}
u|_{\partial\Omega}=0.
\end{equation}

Inside the suspension, chemotaxis is ultimately responsible for the maintenance of the
fluid convection, and thus, for the shape of plumes at large times \cite{CF,LK}.
To propose the the boundary condition for $n$ naturally, we need to consider the effect of chemotaxis on the boundary of suspension.
It's now widely recognized that chemotaxis allows microbial cells to colonize surfaces or interfaces and grow on them,
in the form of multicellular aggregates embedded in matricies commonly referred to as biofilms, which provide the cells with strength in numbers to cope with environmental stresses \cite{Swanson-Reguera-Microbe2016}.
{\it Bacillus subtilis} has long served as a robust model organism to examine the molecular mechanisms of biofilm formation.
Due to the aerotaxis of the cells, {\it Bacillus subtilis} (less domesticated strains) preferably produce biofilms at the air-liquid interface rather than on the surface of a solid phase in a liquid
\cite{Morikawa2006,Swanson-Reguera-Microbe2016}. For the pictures of {\it Bacillus subtilis} biofilms in liquid medium, please refer to Fig.1 in \cite{Morikawa2006} and Fig.2 in \cite{Vlamakis-Chai-Beauregard2013}.
Attachment is initially reversible, and the suspended cell comes and goes until it sticks to the interface and commits to a sessile existence.
However, biofilms are not static entities and cells can be released from the biofilms through an active process of dispersal due to resource limitation and waste product accumulation,
see \cite{Vlamakis-Chai-Beauregard2013,Madigan-Bender-Buckley2019}.

Taking account of the fact that inducing biofilm formation is an important selective advantage of chemotaxis  and biofilm formation is a nearly universal bacterial trait
\cite{Swanson-Reguera-Microbe2016,Vlamakis-Chai-Beauregard2013},
in the present paper, we assume that aerobic bacteria (such as {\it Bacillus subtilis}) may form biofilms at the boundary due to chemotaxis. We also assume that
bacteria may ``cross the boundary'', that is, the cells escape from the system when they are attached at the interface and become nonmotile,
and the cells enter the system when they are released from biofilms.
Therefore, we propose Neumann boundary condition for $n$:
\begin{equation}\label{Neumann boundary condition for n}
\frac{\partial n^m}{\partial \nu}|_{\partial\Omega}=0.
\end{equation}
In view of the cell flux $-\nabla n^m+nu+\chi n\nabla c$, the boundary conditions \eqref{Dirichlet boundary condition for u}--\eqref{Neumann boundary condition for n} imply that there may be three cases occurring on the boundary:
\begin{itemize}
  \item [(1)] If no exchange of oxygen takes place at $\Gamma_1$ on the boundary, that is,
$\frac{\partial c}{\partial \nu}|_{\Gamma_1}=0$ with $\Gamma_1\subset\partial\Omega$, then no bacteria ``cross'' $\Gamma_1$. This case exists at the bottom (and the sides) of the liquid medium
and implies that biofilms can not be formed at these regions.
In particular, this situation includes the no-flux boundary conditions for $c$ and $n$ by letting $\Gamma_1=\partial\Omega$,
which can be considered as an extension of early works; \\[-8mm]
  \item [(2)] If the incoming oxygen molecules are more than leaving oxygen molecules at some region $\Gamma_2\subset\partial\Omega$,
that is, $\frac{\partial c}{\partial \nu}|_{\Gamma_2}>0$, then more bacteria escape from the system at $\Gamma_2$ to produce biofilms;\\[-8mm]
  \item [(3)] If the leaving oxygen molecules are more than incoming oxygen molecules at some region $\Gamma_3\subset\partial\Omega$,
that is, $\frac{\partial c}{\partial \nu}|_{\Gamma_3}<0$, then either more bacteria are released from the biofilms and enter the system,
 or there is no bacteria at the region $\Gamma_3$.
\end{itemize}

Thus, we have the following inhomogeneous mixed boundary conditions
\begin{equation}\label{condition3}
\frac{\partial n^m}{\partial \nu}|_{\partial\Omega}=0, \quad \frac{\partial c}{\partial \nu}|_{\partial\Omega}=\frac{\partial g_1(x,t)}{\partial \nu}(c(x,t)-g_2(x,t)), \quad u|_{\partial\Omega}=0.
 \end{equation}

The purpose of this paper is to establish the existence of bounded weak and strong time periodic solutions
for the problem \eqref{problem0} and \eqref{condition3} in dimension 3.
Since that the boundary condition of oxygen concentration $c$ is inhomogeneous,
it is necessary to make a transformation to apply the standard Neumann heat semigroup argument
and integration by parts to $\eqref{problem0}_2$.
Let
\begin{equation}\label{homogeneous transformation of c}
\tilde{c}=e^{-g_1}(c-g_2).
\end{equation}
Then we have
$$
\frac{\partial \tilde{c}}{\partial \nu}|_{\partial\Omega}=-e^{-g_1}\frac{\partial g_2}{\partial \nu}|_{\partial\Omega}=0.
$$
And the problem \eqref{problem0} and \eqref{condition3} is transformed into
\begin{equation}\label{homogeneous problem}
 \left\{
\begin{array}{l}
n_t+u\cdot\nabla n=\Delta n^m-\chi\nabla\cdot(e^{g_1}n\nabla \tilde{c}+e^{g_1}n\tilde{c}\nabla g_1+n\nabla g_2)+\mu n(a-n)+g, \\[2mm]
\tilde{c}_t-\Delta \tilde{c}+(u-2\nabla g_1)\cdot\nabla \tilde{c}=(|\nabla g_1|^2+\Delta g_1-n-u\nabla g_1-g_{1t})\tilde{c}
\\[1mm]
\hspace{35pt} +(\Delta g_2-u\nabla g_2-ng_2-g_{2t})e^{-g_1}, \\[2mm]
u_t= \Delta u  -\nabla \pi +n\nabla \varphi, \\[2mm]
\nabla\cdot u =0,\\[2mm]
\frac{\partial n^m}{\partial \nu}|_{\partial\Omega}=\frac{\partial \tilde{c}}{\partial \nu}|_{\partial\Omega}=u|_{\partial\Omega}=0.
 \end{array} \right.
\end{equation}
By using a double-level approximation method \cite{Jin2020DCDS}
and some iterative techniques, we obtain the existence and time-space uniform boundedness of
weak time periodic solutions of \eqref{homogeneous problem} for any $m>1$. Moreover,
by deriving the following estimate
$$
\int_{0}^T\int_{\Omega}(n_{\varepsilon}+\varepsilon)^{m-4}|\nabla n_{\varepsilon}|^{4}\,\mathrm{d}x\mathrm{d}s\le C,
$$
where $C$ is independent of $\varepsilon$, we improve the regularity for $m\le\frac{4}{3}$
and show that the obtained periodic solutions are in fact strong periodic solutions.

It is worth mentioning that there exist some essential difficulties to establish the prior estimates of time periodic solutions
for the problem \eqref{homogeneous problem}.

Firstly, considering that  $\eqref{homogeneous problem}_1$ may be degenerate due to $m>1$, and in general does not allow for classical solvability
as the well-known porous medium equations \cite{Vazquez2006}, there may not be enough compactness to get the existences of weak and strong periodic solutions.
In this paper, we use a fourth order regularized problem (see \eqref{fourth order problem} below) to approach the original problem \eqref{homogeneous problem}.
However, different from the second order parabolic system, there is no positivity for the fourth order regularized system.
The most basic and natural $L^1$-norm estimate of $n$ is no longer valid, which brings great difficulties to the later proof. So we
By using a double-level approximation method and introducing three terms $\varepsilon|n|^sn$, $A|n|$ and $An$ in the regularized equation $\eqref{fourth order problem}_1$,
we solve the difficulties caused by the lack of positivity.
For more details, please see the proof of Lemma \ref{solution for linearized problem of n-fourth order} and \ref{priori estimate for fourth order system}.

Secondly, the Neumann boundary condition for $n$ may produce boundary integral at $\partial\Omega$. Indeed,
integrating $\eqref{homogeneous problem}_1$ over $\Omega\times(t_0,t)$ formally,
we see that the following term in the resulting equation
\begin{equation}\label{boundary estimate}
-\chi\int_{t_0}^t\int_{\Omega}\nabla\cdot(e^{g_1}n\tilde{c}\nabla g_1)\,\mathrm{d}x\mathrm{d}s
=-\chi\int_{t_0}^t\int_{\partial\Omega}e^{g_1}n\tilde{c}\frac{\partial g_1}{\partial \nu}\,\mathrm{d}\Gamma\mathrm{d}s\ne0,
\end{equation}
cannot be eliminated or estimated directly, which also prevent us from deriving the positivity, $L^1$-norm and $L^2$-norm estimates of $n$.
By properly choosing the form of chemotactic flux in the regularized equation $\eqref{fourth order problem}_1$, that is,
$$
e^{g_1}n_+\nabla \tilde{c}+e^{g_1}n_+\tilde{c}\nabla g_1+n_+\nabla g_2,
$$
and making some more accurate estimates on \eqref{boundary estimate} by virtue of boundary trace imbedding lemma, we are able to overcome these difficulties.
For more details, please refer to Section 4.

\vspace{3mm}

Now we are the position to give the first result in this paper.

\begin{theorem}\label{main result-1}
Let $m>1$. Assume that \eqref{condition1} and \eqref{condition2} hold, $a, g, \nabla\varphi\in L^{\infty}_T(Q)$ and $g\ge0$.
Then the problem \eqref{homogeneous problem}
admits a bounded weak time periodic solution $(n, c, u)$ with
$n, c\ge0$, $n\in\mathcal{X}_1$, $\tilde{c}\in\mathcal{X}_2$, $u\in\mathcal{X}_3$, where
\begin{align*}
&\mathcal{X}_1=\{n;\ n\in L^{\infty}_T(Q),\ \nabla n^m\in L^{\infty}_T(\mathbb{R}^+, L^2(\Omega)),\ \nabla n^{\frac{m}{2}}\ {\rm and}\  (n^{\frac{m+1}{2}})_t\in L^{2}_T(Q)\},
\\[1mm]
&\mathcal{X}_2=\{\tilde{c};\ \tilde{c}\in L^{\infty}_T(\mathbb{R}^+, W^{1,\infty}(\Omega)),\ \Delta\tilde{c}\ {\rm and}\  \tilde{c}_t\in L^{p}_T(Q) \ {\rm for \ any}\ p>1\},
\\[1mm]
&\mathcal{X}_3=\{u; \ u\in L^{\infty}_T(\mathbb{R}^+, W^{1,\infty}\cap H^1_{\sigma}(\Omega)),\ \Delta u\ {\rm and}\  u_t\in L^{p}_T(Q) \ {\rm for \ any}\ p>1\},
\end{align*}
such that
\begin{align}\label{main result-1-1}
&\sup_{t}(\|n(t)\|_{L^\infty}+\|\tilde{c}(t)\|_{W^{1,\infty}}+\|u(t)\|_{W^{1,\infty}}+\|\nabla n^m(t)\|_{L^2})\le C,
\\\label{main result-1-2}
&\int_0^T(\|\nabla n^{\frac{m}{2}}\|_{L^2}^2+\|(n^{\frac{m+1}{2}})_t\|_{L^2}^2)\,\mathrm{d}s\le C,
\\\label{main result-1-3}
&\int_{0}^T\left(\|u_t\|^p_{L^p}+\|u\|^p_{W^{2,p}}+\|\tilde{c}_t\|^p_{L^p}+\|\tilde{c}\|^p_{W^{2,p}}\right)\,\mathrm{d}s\le C \ \ {\rm for\ any\ } p>1,
\end{align}
where $C$ only depends on $m, \chi, \mu, \Omega, T, a, g, \varphi, p$.
\end{theorem}

\vspace{3mm}

The second result of this paper is concerned with the existence of strong time periodic solutions.

\begin{theorem}\label{main result-2}
Let $1<m\le\frac{4}{3}$. Assume that \eqref{condition1} and \eqref{condition2} hold, $a, g, \nabla\varphi\in L^{\infty}_T(Q)$ and $g\ge0$.
Then the problem \eqref{homogeneous problem}
admits a bounded strong time periodic  solution $(n, c, u)$ with $n, c\ge0$, $n\in\mathcal{D}_1$, $\tilde{c}\in\mathcal{X}_2$, $u\in\mathcal{X}_3$, where
\begin{align*}
&\mathcal{D}_1=\{n;\ n\in L^{\infty}_T(Q),\ \nabla \sqrt{n}\in L^{\infty}_T(\mathbb{R}^+, L^2(\Omega)),\ \Delta n^m, \ \nabla n^{\frac{m}{4}}\ {\rm and}\  n_t\in L^{2}_T(Q)\},
\end{align*}
such that
\begin{align}\label{main result-2-1}
&\sup_{t}(\|n(t)\|_{L^\infty}+\|\tilde{c}(t)\|_{W^{1,\infty}}+\|u(t)\|_{W^{1,\infty}}+\|\nabla \sqrt{n}(t)\|_{L^2})\le C,
\\\label{main result-2-2}
&\int_0^T(\|\Delta n^m\|_{L^2}^2+\|\nabla n^{\frac{m}{4}}\|_{L^2}^2+\|n_t\|_{L^2}^2)\,\mathrm{d}s\le C,
\\\label{main result-2-3}
&\int_{0}^T\left(\|u_t\|^p_{L^p}+\|u\|^p_{W^{2,p}}+\|\tilde{c}_t\|^p_{L^p}+\|\tilde{c}\|^p_{W^{2,p}}\right)\,\mathrm{d}s\le C \ \ {\rm for\ any\ } p>1,
\end{align}
where $C$ only depends on $m, \chi, \mu, \Omega, T, a, g, \varphi, p$.
\end{theorem}

The remainder of this paper is organized as follows. In Section 2,
we recall some auxiliary lemmas which will be used in this paper.
In Section 3, we prove the existence of time periodic solutions for a
fourth order regularized problem which approaches the original problem \eqref{homogeneous problem}.
In section 4, by using a double-level approximation method and some iterative techniques,
we obtain the existence and time-space uniform boundedness of
weak time periodic solutions for any $m>1$.
Finally, we show that the obtained time periodic solutions are in fact strong periodic solutions
for $m\le\frac{4}{3}$ in the last section.

\vspace{3mm}

\setcounter{equation}{0}
\section{Some auxiliary lemmas}

{\bf Notations}.
\begin{itemize}
  \item $Q=\Omega\times\mathbb{R}^+$ and $Q_T=\Omega\times(0,T)$.\\[-7.5mm]
  \item $f\in L^{p}_T(\mathbb{R}^+; X)\Leftrightarrow f$ is a time periodic function with period $T$, and $f\in L^{p}(0,T; X)$.
        For simplicity, we denote $L^{p}_T(\mathbb{R}^+; L^{p}(\Omega))$ by $L^{p}_T(Q)$.\\[-7.5mm]
  \item The outward unit normal to $\partial\Omega$ is denoted by $\nu$.\\[-7.5mm]
  \item $C^{\infty}_{0,\sigma}(\Omega)$ denotes the set of all $C^{\infty}_{0,\sigma}(\Omega)$-real functions $\varphi=(\varphi_1, \varphi_2, \varphi_3)$ with compact support in $\Omega$, such that
$\nabla\cdot \varphi =0$. The closure of $C^{\infty}_{0,\sigma}(\Omega)$ with respect to norm $L^r$ is denoted by $L^r_{\sigma}(\Omega)$.\\[-7.5mm]
  \item $C$ stands for a generic positive constant which may vary from line to line.
\end{itemize}

By \cite{Galdi1994}, each $u\in L^r(\Omega)$ has a unique decomposition
$$
u=v+\nabla p, \quad v\in L^r_{\sigma},\quad \nabla p\in G^r
$$
with $G^r=\{\nabla p; \nabla p\in L^r, p\in L^r_{loc}\}$, and the projection
$P: L^r(\Omega)\to L^r_{\sigma}(\Omega)$ is called Helmholtz projection.
Let $\mathcal{A}\omega:=-P\Delta\omega$, then $\mathcal{A}$ generates a bounded analytic semigroup $\{e^{-t\mathcal{A}}; t\ge0\}$ on $L^r_{\sigma}(\Omega)$,
and the time periodic solution $u$ of \eqref{problem0} can be expressed as
\begin{equation}\label{Navier-Stokes}
u=\int_{-\infty}^te^{-(t-s)\mathcal{A}}P(n\nabla\varphi)\,\mathrm{d}s.
\end{equation}
For more details, please refer to \cite{Farwig-Okabe2010,Kozono-Yanagisawa2009}.

By \cite{Jin2020PRSE}, we have the following two lemmas.

\begin{lemma}\label{periodic gronwall inequality-1}
Let $T > 0, a > 0, \sigma\ge0$, and suppose that $f: \mathbb{R}^+\to[0,\infty)$ is absolutely
continuous, $f, h$ are time periodic functions with period $T$, and $f$ satisfies
\begin{equation*}
f(t)-f(t_0)+a\int_{t_0}^tf^{1+\sigma}(s)\,\mathrm{d}s\le\int_{t_0}^th(s)\,\mathrm{d}s,
\quad 0\le t_0< t,
\end{equation*}
where $0\le f, h\in L^1_T(\mathbb{R}^+)$ and $\int_0^Th(s)\,\mathrm{d}s\le \beta$.
Then we have
$$
\sup_{t\in(0,T)}f(t)+a\int_0^Tf(t)\,\mathrm{d}t\le (\frac{\beta}{aT})^{1/(1+\sigma)}+2\beta.
$$
\end{lemma}

\begin{lemma}\label{periodic gronwall inequality-2}
Let $T > 0, a > 0, \sigma>0$, and suppose that $f: \mathbb{R}^+\to[0,\infty)$ is absolutely
continuous, $f, g, h$ are time periodic functions with period $T$, and satisfy
\begin{equation*}
f(t)-f(t_0)+a\int_{t_0}^tf^{1+\sigma}(s)\,\mathrm{d}s\le\int_{t_0}^tg(s)f(s)\,\mathrm{d}s+\int_{t_0}^th(s)\,\mathrm{d}s,
\quad 0\le t_0< t,
\end{equation*}
where $0\le g, h\in L^1_T(\mathbb{R}^+)$, $\int_0^Tg(s)\,\mathrm{d}s\le \alpha$ and $\int_0^Th(s)\,\mathrm{d}s\le \beta$.
Then we have
$$
\sup_{t\in(0,T)}f(t)+a\int_0^Tf^{1+\sigma}(t)\,\mathrm{d}t\le C,
$$
where $C$ is a constant depending only on $a, \alpha, \beta, T$. While, if $a = 0$ and $\int_0^Tf(s)\,\mathrm{d}s\le \gamma$.
Then we also have
$$
\sup_{t\in(0,T)}f(t)\le C,
$$
where $C$ is a constant depending only on $\gamma, \alpha, \beta, T$.
\end{lemma}

By \cite{Jin2017ZAMP,Nakao-Koyanagi1985}, we also have the following lemma.

\begin{lemma}\label{L^pL^q}
Assume that $f\in L^{p}_T(\mathbb{R}^+; L^p(\Omega))$ with $p>1$. Then the following problem
\begin{equation}\label{problem1}
 \left\{
\begin{array}{l}
u_t-\Delta u+u=f(x,t), \\[1.5mm]
\frac{\partial u}{\partial \nu}|_{\partial\Omega}=0
 \end{array} \right.
 \end{equation}
admits a unique strong time periodic solution $u\in W^{2,1}_2(Q_T)$, and
$$
\|u\|^p_{W^{2,1}_{p}(Q_T)} \le C\|f\|_{L^p(Q_T)}^p,
$$
where $C$ is a positive constant.
\end{lemma}

By Gagliardo-Nirenberg interpolation inequality\cite{Leoni2017}, we see that

\begin{lemma}\label{G-N}
Let $m, k$ be non-negative integers, $1\le p, q, r, s\le \infty$ and $\Omega$ be a bounded Lipschitz domain in $\mathbb{R}^N$.
Suppose that $\frac{1}{q}$ is between $\frac{1}{p}-\frac{m-k}{N}$ and $\frac{k}{mp}+\frac{m-k}{mr}$, and when $\frac{1}{q}=\frac{1}{p}-\frac{m-k}{N}$,
$m-k-\frac{N}{p}$ is not non-negative integer. Then if $u\in L^r(\Omega)\cap L^s(\Omega)$ and $\partial^m u\in L^p(\Omega)$, we have $\partial^ku\in L^q(\Omega)$ and
\begin{equation}\label{G-N-1}
\|\partial^ku\|_{L^q(\Omega)}\le C\|u\|^{\theta}_{L^r(\Omega)}\|\partial^m u\|^{1-\theta}_{L^p(\Omega)}+C\|u\|_{L^s(\Omega)},
\end{equation}
where $C>0$ depends only on $N, m, k, p, q, r, s, \Omega$, and $\theta\in[0,\frac{m-k}{m}]$ satisfying
$$
\frac{1}{q}=\theta(\frac{1}{r}+\frac{k}{N})+(1-\theta)(\frac{1}{p}-\frac{m-k}{N}).
$$
Especially, if $u|_{\partial\Omega}=0$, then \eqref{G-N-1} is reduced to
\begin{equation}\label{G-N-2}
\|\partial^ku\|_{L^q(\Omega)}\le C\|u\|^{\theta}_{L^r(\Omega)}\|\partial^m u\|^{1-\theta}_{L^p(\Omega)}.
\end{equation}
\end{lemma}

The following boundary trace imbedding lemma follows from \cite{Adams-Fournier}.

\begin{lemma}\label{Trace Imbedding}
Let $m$ be a non-negative integer, $1\le p, q<\infty$ and $\Omega$ be a domain in $\mathbb{R}^N$ satisfying the uniform $C^m$-regularity condition, and suppose there exists
a simple $(m,p)$-extension operator $E$ for $\Omega$. Also suppose that $mp<N$ and $p\le q\le \frac{(N-1)p}{N-mp}$. Then
\begin{equation}\label{Trace Imbedding-1}
W^{m,p}(\Omega)\hookrightarrow L^q(\partial\Omega).
\end{equation}
If $mp=N$, then imbedding \eqref{Trace Imbedding-1} holds for $1\le p\le q<\infty$.
\end{lemma}

\begin{lemma}[\cite{Mizoguchi-Souplet2014}]\label{outward normal derivative}
Assume that $\Omega$ is bounded and $\omega\in C^2(\overline{\Omega})$ satisfying $\frac{\partial \omega}{\partial \nu}|_{\partial\Omega}=0$.
Then we get that
$$
\frac{\partial |\nabla\omega|^2}{\partial \nu}\le 2\kappa|\nabla\omega|^2 \ \ {\rm on}\ \ \partial\Omega,
$$
where $\kappa>0$ is an upper bound for the curvature of $\Omega$.
\end{lemma}

\vspace{5mm}

\setcounter{equation}{0}
\section{Time periodic solutions for a fourth-order regularized problem}

To obtain the compactness of the operator, we
use a fourth order regularized system as follows to approach the original system.
\begin{equation}\label{fourth order problem}
 \left\{
\begin{array}{l}
n_t-m\nabla\cdot((|n|+\varepsilon)^{m-1}\nabla n)+\delta\Delta^2n+\varepsilon|n|^sn+u\cdot\nabla n+An
\\[1mm]
\hspace{35pt} =-\chi\nabla\cdot(e^{g_1}n_+\nabla \tilde{c}+e^{g_1}n_+\tilde{c}\nabla g_1+n_+\nabla g_2)+\mu |n|(a-n)+A|n|+g, \\[2mm]
\tilde{c}_t-\Delta \tilde{c}+(u-2\nabla g_1)\cdot\nabla \tilde{c}=(|\nabla g_1|^2+\Delta g_1-n_+-u\nabla g_1-g_{1t})\tilde{c}
\\[1mm]
\hspace{35pt} +(\Delta g_2-u\nabla g_2-ng_2-g_{2t})e^{-g_1}, \\[2mm]
u_t= \Delta u  -\nabla \pi +n\nabla \varphi, \\[2mm]
\nabla\cdot u =0,\\[2mm]
\frac{\partial n}{\partial \nu}|_{\partial\Omega}=\frac{\partial \Delta n}{\partial \nu}|_{\partial\Omega}
=\frac{\partial \tilde{c}}{\partial \nu}|_{\partial\Omega}=u|_{\partial\Omega}=0,
 \end{array} \right.
 \end{equation}
where $\max\{2(m-1), 2\}<s\le5m-1$, $\delta, \varepsilon>0$, $m>1$, $A>0$. The terms $\varepsilon|n|^sn$, $An$ and $A|n|$ are introduced in
$\eqref{fourth order problem}_1$ to solve the difficulties caused by the lack of positivity
of the fourth order regularized problem.

In this section, we prove the existence of time periodic solutions for the fourth-order problem \eqref{fourth order problem}.
For this purpose, we linearize this problem.
Consider the following problem
\begin{equation}\label{linearized problem of u}
 \left\{
\begin{array}{l}
u_t-\Delta u +\nabla \pi =\eta\hat{n}\nabla \varphi, \\[2mm]
\nabla\cdot u =0,\\[2mm]
u|_{\partial\Omega}=0,
 \end{array} \right.
 \end{equation}
where $\eta\in[0,1]$ is a constant. By \cite{Jin2017ZAMP}, we have

\begin{lemma}\label{solution for linearized problem of u}
Assume $\nabla \varphi\in L^{\infty}_T(Q)$, $\hat{u}\in L^4_T(\mathbb{R}^+; L^4_{\sigma}(\Omega))$ and $\hat{n}\in L^2_T(\mathbb{R}^+; L^2(\Omega))$.
Then \eqref{linearized problem of u}
admits a unique strong time periodic solution $u$ with $u\in L^{\infty}(0,T;H^1_{\sigma}(\Omega))\cap L^{2}(0,T;H^2_{\sigma}(\Omega))$
and $u_t\in L^{2}(0,T;L^2_{\sigma}(\Omega))$.
\end{lemma}

For the above solution $u$, we consider the following problem for any $\eta\in[0,1]$.
\begin{equation}\label{linearized problem of c}
 \left\{
\begin{array}{l}
c_t-\Delta c+u\cdot\nabla c+(1-\eta)c=-\hat{n}_+c, \\[2mm]
\frac{\partial c}{\partial \nu}|_{\partial\Omega}=\eta\frac{\partial g_1(x,t)}{\partial \nu}(c(x,t)-g_2(x,t)).
 \end{array} \right.
 \end{equation}
Let $\tilde{c}=e^{-\eta g_1}(c-g_2)$.
Then \eqref{linearized problem of c} is equivalent to
\begin{equation}\label{linearized problem of tilde-c}
 \left\{
\begin{array}{l}
\tilde{c}_t-\Delta \tilde{c}+(u-2\eta\nabla g_1)\cdot\nabla \tilde{c}+(1-\eta+\hat{n}_+)\tilde{c}=\eta(\eta|\nabla g_1|^2+\Delta g_1-u\nabla g_1-g_{1t})\tilde{c}
\\[1mm]
\hspace{35pt} +(\Delta g_2-u\nabla g_2-\hat{n}_+g_2-g_{2t}-(1-\eta)g_2)e^{-\eta g_1}, \\[2mm]
\frac{\partial \tilde{c}}{\partial \nu}|_{\partial\Omega}=0.
 \end{array} \right.
 \end{equation}

By \cite{Jin2017ZAMP}, we have
\begin{lemma}\label{solution for linearized problem of tilde-c}
Assume $\hat{n}\in L^2_T(\mathbb{R}^+; H^1(\Omega))$.
Let $u$ be the time periodic solution of the problem \eqref{linearized problem of u}.
Then \eqref{linearized problem of tilde-c} (or \eqref{linearized problem of c}) admits a unique strong time periodic solution $\tilde{c}$
with $c\ge0$, $\tilde{c}\in L^{\infty}_T(Q)\cap L^{\infty}_T(\mathbb{R}^+;H^2(\Omega))\cap L^{2}_T(\mathbb{R}^+;H^3(\Omega))$
and $\tilde{c}_t\in L^{\infty}_T(\mathbb{R}^+;L^2(\Omega))$.
\end{lemma}

For the above obtained solutions $u, \tilde{c}$, we consider the following problem.
\begin{equation}\label{linearized problem of n-fourth order}
 \left\{
\begin{array}{l}
n_t-m\nabla\cdot((|\hat{n}|+\varepsilon)^{m-1}\nabla n)+\delta\Delta^2n+\varepsilon|\hat{n}|^sn+An+u\cdot\nabla n
\\[1mm]
=-\eta\chi\nabla\cdot(e^{g_1}n_+\nabla \tilde{c}+e^{g_1}n_+\tilde{c}\nabla g_1+n_+\nabla g_2)+\eta(\mu a+A)|\hat{n}|-\mu|\hat{n}|n+\eta g, \\[2mm]
\frac{\partial n}{\partial \nu}|_{\partial\Omega}=\frac{\partial \Delta n}{\partial \nu}|_{\partial\Omega}=0.
 \end{array} \right.
 \end{equation}

For the above linear parabolic problem, when $A$ is sufficiently large, the existence
of time periodic solutions can be easily obtained by a fixed point method. That
is, define a Poincar\'{e} map from $n(x, 0)$ to $n(x, T)$, the time-periodic solution is
then identified as a fixed point of this Poincar\'{e} map (see also \cite{Jin2020DCDS}). We only give the regularity
estimates.

For simplicity, in what follows, we may assume that the solution $n$ is sufficiently
smooth, otherwise, we can approximate $u, \tilde{c}, n$ with a sequence of sufficiently
smooth functions $u_k, \tilde{c}_k, n_k$ such that the corresponding solutions n are sufficiently
smooth, and the following energy estimates can be obtained through an
approximate process.

\begin{lemma}\label{solution for linearized problem of n-fourth order}
Assume $a, g\in L^{\infty}_T(Q)$, $\hat{n}\in L^{\infty}_T(Q)\cap L^{\infty}_T(\mathbb{R}^+; H^1(\Omega))\cap L^2_T(\mathbb{R}^+; H^2(\Omega))$
and $\frac{\partial \hat{n}}{\partial \nu}|_{\partial\Omega}=0$.
Let $u, \tilde{c}, n$ be the time periodic solutions of the problem \eqref{linearized problem of u}, \eqref{linearized problem of tilde-c}
and \eqref{linearized problem of n-fourth order}, respectively. If $A>1$ is a sufficiently large constant,
then $n\in L^{\infty}_T(\mathbb{R}^+;H^2(\Omega))\cap L^{2}_T(\mathbb{R}^+;H^4(\Omega))$
and $n_t\in L^{2}_T(\mathbb{R}^+;L^2(\Omega))$.
\end{lemma}

\noindent{\bf Proof}.\
Multiplying $\eqref{linearized problem of n-fourth order}_1$ by $n$, integrating it over $\Omega\times(t_0,t)$ for
any $t_0< t\le t_0+T$, and using Lemma \ref{G-N}, \ref{Trace Imbedding} and \ref{solution for linearized problem of tilde-c},
we see that for sufficiently large $A$,
\begin{align*}
&\frac{1}{2}\int_{\Omega}(|n(x,t)|^2-|n(x,t_0)|^2)\,\mathrm{d}x+\int_{t_0}^t\int_{\Omega}(\delta|\Delta n|^2+\varepsilon|\hat{n}|^sn^2+\mu|\hat{n}|n^2+An^2)\,\mathrm{d}x\mathrm{d}s
   \nonumber
\\
& \le  \eta\chi\int_{t_0}^t\int_{\Omega}(e^{g_1}n_+\nabla \tilde{c}+e^{g_1}n_+\tilde{c}\nabla g_1+n_+\nabla g_2)\nabla n\,\mathrm{d}x\mathrm{d}s
-\eta\chi\int_{t_0}^t\int_{\partial\Omega}e^{g_1}n_+n\tilde{c}\frac{\partial g_1}{\partial \nu}\,\mathrm{d}\Gamma\mathrm{d}s  \nonumber
\\
& \hspace{12pt} +\int_{t_0}^t\int_{\Omega}(\eta(\mu a+A)|\hat{n}|-\mu|\hat{n}|n+\eta g)n\,\mathrm{d}x\mathrm{d}s  \nonumber
\\
& \le  C\int_{t_0}^t\|n\|_{L^6}\|\nabla\tilde{c}\|_{L^3}\|\nabla n\|_{L^2}\,\mathrm{d}s+ C\int_{t_0}^t\|n\|_{L^2}\|\nabla n\|_{L^2}\,\mathrm{d}s
+C\int_{t_0}^t\|n\|^2_{L^2(\partial\Omega)}\,\mathrm{d}s   \nonumber
\\
&\hspace{12pt} +\int_{t_0}^t\|(\eta(\mu a+A)|\hat{n}|-\mu|\hat{n}|n+\eta g)\|_{L^2}\|n\|_{L^2}\,\mathrm{d}s  \nonumber
\\
& \le  C\int_{t_0}^t\|\nabla n\|^2_{L^2}\,\mathrm{d}s+ \frac{A}{2}\int_{t_0}^t\|n\|^2_{L^2}\,\mathrm{d}s+C,
\end{align*}
which implies
\begin{align*}
\int_{\Omega}(|n(x,t)|^2-|n(x,t_0)|^2)\,\mathrm{d}x+\int_{t_0}^t\int_{\Omega}(2\delta|\Delta n|^2+An^2)\,\mathrm{d}x\mathrm{d}s
\le  C\int_{t_0}^t\|\nabla n\|^2_{L^2}\,\mathrm{d}s+C.
\end{align*}
By Lemma \ref{periodic gronwall inequality-1}, it follows that
\begin{align}\label{solution for linearized problem of n-fourth order-1}
\sup_{t}\|n(t)\|^2_{L^2}+\int_{0}^T\|n(s)\|^2_{H^2}\,\mathrm{d}s
\le  C\int_{0}^T\|\nabla n\|^2_{L^2}\,\mathrm{d}s+C,
\end{align}
where the constant $C>0$ is independent of $A$.
Multiplying $\eqref{linearized problem of n-fourth order}_1$ by $(-\Delta)n$ and integrating it over $\Omega\times(t_0,t)$,
and using Lemma \ref{G-N} and \ref{solution for linearized problem of tilde-c},
we have
\begin{align*}
&\frac{1}{2}\int_{\Omega}(|\nabla n(x,t)|^2-|\nabla n(x,t_0)|^2)\,\mathrm{d}x+\int_{t_0}^t\int_{\Omega}(\delta|\nabla\Delta n|^2+A|\nabla n|^2)\,\mathrm{d}x\mathrm{d}s
   \nonumber
\\
& =   \int_{t_0}^t\int_{\Omega}\left(m(|\hat{n}|+\varepsilon)^{m-1}\nabla n-\eta\chi e^{g_1}n_+\nabla \tilde{c}-\eta\chi n_+\nabla g_2\right)\cdot\nabla\Delta n\,\mathrm{d}x\mathrm{d}s   \nonumber
\\
& \hspace{12pt} +\eta\chi\int_{t_0}^t\int_{\Omega}\nabla\cdot(e^{g_1}n_+\tilde{c}\nabla g_1)\Delta n\,\mathrm{d}x\mathrm{d}s
+\int_{t_0}^t\int_{\Omega}u\cdot\nabla n\Delta n\,\mathrm{d}x\mathrm{d}s  \nonumber
\\
& \hspace{12pt} +\int_{t_0}^t\int_{\Omega}(\varepsilon|\hat{n}|^sn+\mu|\hat{n}|n-\eta(\mu a+A)|\hat{n}|-\eta g)\Delta n\,\mathrm{d}x\mathrm{d}s  \nonumber
\\
& \le  C\int_{t_0}^t\left(\|\nabla n\|_{L^2}+\|n\|_{L^{\infty}}\|\nabla\tilde{c}\|_{L^2}+\| n\|_{L^2}\right)\|\nabla\Delta n\|_{L^2}\,\mathrm{d}s   \nonumber
\\
&\hspace{12pt} +\eta\chi\int_{t_0}^t\|\nabla\cdot(e^{g_1}n_+\tilde{c}\nabla g_1)\|_{L^2}\|\Delta n\|_{L^2}\,\mathrm{d}s+\int_{t_0}^t\|u\|_{L^6}\|\nabla n\|_{L^3}\|\Delta n\|_{L^2}\,\mathrm{d}s  \nonumber
\\
&\hspace{12pt} +\int_{t_0}^t\|(\varepsilon|\hat{n}|^sn+\mu|\hat{n}|n-\eta(\mu a+A)|\hat{n}|-\eta g)\|_{L^2}\|\Delta n\|_{L^2}\,\mathrm{d}s  \nonumber
\\
& \le  C\int_{t_0}^t\left(\|n\|_{H^1}+\|n\|_{H^2}\|\nabla\tilde{c}\|_{L^2}\right)\|\nabla\Delta n\|_{L^2}\,\mathrm{d}s+ C\int_{t_0}^t\|n\|_{H^1}\|\Delta n\|_{L^2}(1+\|\tilde{c}\|_{H^1})\,\mathrm{d}s   \nonumber
\\
&\hspace{12pt} +C\int_{t_0}^t\|u\|_{H^1}\|\nabla n\|_{H^1}\|\Delta n\|_{L^2}\,\mathrm{d}s+C\int_{t_0}^t(\|n\|_{L^2}+1)\|\Delta n\|_{L^2}\,\mathrm{d}s  \nonumber
\\
& \le  \frac{\delta}{2}\int_{t_0}^t\|\nabla\Delta n\|^2_{L^2}\,\mathrm{d}s+ \frac{A}{2}\int_{t_0}^t\|\nabla n\|^2_{L^2}\,\mathrm{d}s
+C\int_{t_0}^t(\|n\|^2_{L^2}+\|\Delta n\|^2_{L^2})\,\mathrm{d}s+C.
\end{align*}
Then,
\begin{align*}
&\int_{\Omega}(|\nabla n(x,t)|^2-|\nabla n(x,t_0)|^2)\,\mathrm{d}x+\int_{t_0}^t\int_{\Omega}(\delta|\nabla\Delta n|^2+A|\nabla n|^2)\,\mathrm{d}x\mathrm{d}s  \nonumber
\\
&\le  C\int_{t_0}^t(\|n\|^2_{L^2}+\|\Delta n\|^2_{L^2})\,\mathrm{d}s+C.
\end{align*}
If $A$ is large enough, combining Lemma \ref{periodic gronwall inequality-1} with inequality \eqref{solution for linearized problem of n-fourth order-1}, we obtain
\begin{align}\label{solution for linearized problem of n-fourth order-2}
\sup_{t}\|\nabla n(t)\|^2_{L^2}+\int_{0}^T\|n(s)\|^2_{H^3}\,\mathrm{d}s
\le  C,
\end{align}
where the constant $C$ depends on $\delta$ and $\varepsilon$.
Multiplying $\eqref{linearized problem of n-fourth order}_1$ by $\Delta^2n$ and integrating it over $\Omega\times(t_0,t)$,
and using the estimate \eqref{solution for linearized problem of n-fourth order-2}, Lemma \ref{G-N} and \ref{solution for linearized problem of tilde-c},
we also have
\begin{align*}
&\frac{1}{2}\int_{\Omega}(|\Delta n(x,t)|^2-|\Delta n(x,t_0)|^2)\,\mathrm{d}x
+\int_{t_0}^t\int_{\Omega}(\delta|\Delta^2 n|^2+A|\Delta n|^2)\,\mathrm{d}x\mathrm{d}s
   \nonumber
\\
& =   \int_{t_0}^t\int_{\Omega}\left(m\nabla\cdot((|\hat{n}|+\varepsilon)^{m-1}\nabla n)-u\cdot\nabla n\right)\Delta^2 n\,\mathrm{d}x\mathrm{d}s   \nonumber
\\
& \hspace{12pt} -\eta\chi\int_{t_0}^t\int_{\Omega}\nabla\cdot(e^{g_1}n_+\nabla \tilde{c}+e^{g_1}n_+\tilde{c}\nabla g_1+n_+\nabla g_2)\Delta^2 n\,\mathrm{d}x\mathrm{d}s  \nonumber
\\
& \hspace{12pt} -\int_{t_0}^t\int_{\Omega}(\varepsilon|\hat{n}|^sn+\mu|\hat{n}|n-\eta(\mu a+A)|\hat{n}|-\eta g)\Delta^2 n\,\mathrm{d}x\mathrm{d}s  \nonumber
\\
& \le  C\int_{t_0}^t\left(\|\Delta n\|_{L^2}+\|\nabla\hat{n}\cdot\nabla n\|_{L^2}+\|u\|_{L^6}\|\nabla n\|_{L^3}\right)\|\Delta^2 n\|_{L^2}\,\mathrm{d}s   \nonumber
\\
&\hspace{12pt} +C\int_{t_0}^t\left(\|n\nabla\tilde{c}\|_{L^2}+\|n\Delta\tilde{c}\|_{L^2}+\|\nabla n\cdot\nabla\tilde{c}\|_{L^2}+\|n\|_{L^2}+1\right)\|\Delta^2 n\|_{L^2}\,\mathrm{d}s  \nonumber
\\
& \le  \frac{\delta}{2}\int_{t_0}^t\|\Delta^2 n\|^2_{L^2}\,\mathrm{d}s+C\int_{t_0}^t\left(\|\Delta n\|^2_{L^2}+\|\nabla\hat{n}\|_{L^4}^2\|\nabla n\|^2_{L^4}+\|u\|^2_{H^1}\|\nabla n\|^2_{H^1}\right)\,\mathrm{d}s   \nonumber
\\
&\hspace{12pt} +C\int_{t_0}^t\left(\|n\|^2_{L^\infty}\|\tilde{c}\|^2_{L^2}+\|n\|^2_{L^{\infty}}\|\Delta\tilde{c}\|^2_{L^2}+\|\nabla n\|_{L^4}^2\|\nabla\tilde{c}\|^2_{L^4}+1\right)\,\mathrm{d}s  \nonumber
\\
& \le  \frac{\delta}{2}\int_{t_0}^t\|\Delta^2 n\|^2_{L^2}\,\mathrm{d}s+C(\sup_t\|u\|^2_{H^1}+\sup_t\|\tilde{c}\|^2_{H^2}+1)\int_{t_0}^t\|n\|^2_{H^2}\,\mathrm{d}s+C.
\end{align*}
Combining with Lemma \ref{periodic gronwall inequality-1}, we finally get
\begin{align}\label{solution for linearized problem of n-fourth order-3}
\sup_{t}\|n(t)\|^2_{H^2}+\int_{0}^T\|n(s)\|^2_{H^4}\,\mathrm{d}s
\le  C,
\end{align}
where the constant $C$ depends on $\delta$ and $\varepsilon$.
Similarly, multiplying $\eqref{linearized problem of n-fourth order}_1$ by $n_t$ and integrating it, we can obtain
\begin{align}\label{solution for linearized problem of n-fourth order-4}
\int_{0}^T\|n_t(s)\|^2_{L^2}\,\mathrm{d}s
\le  C.
\end{align}
Summing up, we complete the proof.
\hfill$\Box$

\vbox to 3mm{}

Next, we show the existence of time periodic solutions in
dimension 3 by Leray-Schauder's fixed point theorem. Define an operator $\mathcal{F}: \mathcal{E}\times[0,1]\to\mathcal{E}$ as follows:
$$
\mathcal{F}(\hat{n}, \eta)=n,
$$
where
$$
\mathcal{E}=\{n; n\in L^{\infty}_T(Q)\cap L^{\infty}_T(\mathbb{R}^+; H^1(\Omega))\cap L^2_T(\mathbb{R}^+; H^2(\Omega))\},
$$
endowed with the norm
$$
\|n\|_{\mathcal{E}}=\sup_t(\|u(\cdot,t)\|_{L^\infty}+\|n(\cdot, t)\|_{H^1})+\int_0^T\|n(\cdot, t)\|^2_{H^2}\,\mathrm{d}t
$$
and $u, \tilde{c}, n$ are the time periodic solutions of the problem \eqref{linearized problem of u}, \eqref{linearized problem of tilde-c}
and \eqref{linearized problem of n-fourth order}, respectively.

\begin{lemma}\label{priori estimate for fourth order system}
Assume $s>\max\{2m-2,2\}$, $a, g, \nabla\varphi\in L^{\infty}_T(Q)$ and let $\mathcal{F}(n, \eta)=n$ with $\eta\in(0,1]$. Then there exists $R>0$ such that
$$
\|n\|_{\mathcal{E}}\le R,
$$
where $R$ depends on $\delta, \varepsilon$, and is independent of $A$.
\end{lemma}

\noindent{\bf Proof}.\
Take $\hat{n}=n$ in \eqref{linearized problem of n-fourth order}.
Multiplying $\eqref{linearized problem of n-fourth order}_1$ by $n$, integrating it over $\Omega\times(t_0,t)$ for
any $t_0< t\le t_0+T$, and using Lemma \ref{G-N}, \ref{Trace Imbedding} and \ref{solution for linearized problem of tilde-c}, we see that
\begin{align*}
&\frac{1}{2}\int_{\Omega}(|n(x,t)|^2-|n(x,t_0)|^2)\,\mathrm{d}x+\int_{t_0}^t\int_{\Omega}m(|n|+\varepsilon)^{m-1}|\nabla n|^2\,\mathrm{d}x\mathrm{d}s
\\
& \hspace{12pt} +\int_{t_0}^t\int_{\Omega}(\delta|\Delta n|^2+\varepsilon|n|^{s+2}+\mu|n|^3+An^2)\,\mathrm{d}x\mathrm{d}s
   \nonumber
\\
& \le  -\frac{\eta\chi}{2}\int_{t_0}^t\int_{\Omega}(\nabla e^{g_1}\nabla \tilde{c}+e^{g_1}\Delta\tilde{c}+\Delta g_2)n_+^2\,\mathrm{d}x\mathrm{d}s
-\eta\chi\int_{t_0}^t\int_{\Omega}\nabla\cdot(e^{g_1}n_+\tilde{c}\nabla g_1)n\,\mathrm{d}x\mathrm{d}s  \nonumber
\\
& \hspace{12pt} +\int_{t_0}^t\int_{\Omega}\eta((\mu a+A)|n|+g)n\,\mathrm{d}x\mathrm{d}s  \nonumber
\\
& \le  C\eta\int_{t_0}^t(\|n\|^4_{L^4}+\|\nabla\tilde{c}\|_{H^1}^2)\,\mathrm{d}s+\frac{m\varepsilon^{m-1}}{2}\int_{t_0}^t\|\nabla n\|^2_{L^2}\,\mathrm{d}s+A\eta\int_{t_0}^t\|n\|^2_{L^2}\,\mathrm{d}s
\\
& \hspace{12pt} +\frac{\mu}{2}\int_{t_0}^t\|n\|^3_{L^3}\,\mathrm{d}s
+C\eta,
\end{align*}
which implies
\begin{align}\label{priori estimate for fourth order system-1}
&\int_{\Omega}(|n(x,t)|^2-|n(x,t_0)|^2)\,\mathrm{d}x+\int_{t_0}^t\int_{\Omega}m(|n|+\varepsilon)^{m-1}|\nabla n|^2\,\mathrm{d}x\mathrm{d}s  \nonumber
\\
& \hspace{12pt} +\int_{t_0}^t\int_{\Omega}(\delta|\Delta n|^2+\varepsilon|n|^{s+2}+\mu|n|^3)\,\mathrm{d}x\mathrm{d}s  \nonumber
\\
&\le  C\eta\int_{t_0}^t(\|n\|^4_{L^4}+\|\nabla\tilde{c}\|_{H^1}^2)\,\mathrm{d}s
+C\eta.
\end{align}
Taking $\hat{n}=n$ in \eqref{linearized problem of u} and multiplying $\eqref{linearized problem of u}_1$ by $u$ and $u_t$,
respectively, then combining the two resulting inequalities and applying Lemma \ref{periodic gronwall inequality-1}, we have
\begin{align}\label{3-11}
\sup_{t}\|u(t)\|^2_{H^1}+\int_{0}^T\left(\|u(s)\|^2_{H^1}+\|u_t(s)\|^2_{L^2}\right)\mathrm{d}s\,\mathrm{d}s
\le  C\int_{0}^T\|n(s)\|^2_{L^2}\,\mathrm{d}s.
\end{align}
Noticing that
$$
-\Delta u+\nabla\pi=-u_t+\eta n\nabla\varphi, \qquad  u|_{\partial\Omega}=0,
$$
using $L^2$-theory of Stokes operator, we have
$$
\|u\|^2_{H^2}\le C (\|u_t\|^2_{L^2}+\|n\|^2_{L^2}),
$$
that is
$$
\int_0^T\|u(s)\|^2_{H^2}\mathrm{d}s\le C \int_0^T(\|u_t(s)\|^2_{L^2}+\|n(s)\|^2_{L^2})\mathrm{d}s.
$$
Combining with \eqref{3-11} yields
\begin{align}\label{priori estimate for fourth order system-2}
\sup_{t}\|u(t)\|^2_{H^1}+\int_{0}^T\left(\|u(s)\|^2_{H^2}+\|u_t(s)\|^2_{L^2}\right)ds\,\mathrm{d}s
\le  C\int_{0}^T\|n\|^2_{L^2}\,\mathrm{d}s.
\end{align}
Take $\hat{n}=n$ in \eqref{linearized problem of tilde-c}. We first multiply $\eqref{linearized problem of tilde-c}_1$ by $\tilde{c}$,
integrate it over $\Omega\times(t_0,t)$ for any $t_0< t\le t_0+T$,
and apply Lemma \ref{periodic gronwall inequality-1} to the resulting inequality. It is easy to see that
\begin{align}\label{priori estimate for fourth order system-3}
\sup_{t}\|\tilde{c}(t)\|^2_{L^2}+\int_{0}^T\|\nabla \tilde{c}(s)\|^2_{L^2}\,\mathrm{d}s
\le  C\int_{0}^T\|n(s)\|^2_{L^2}\,\mathrm{d}s.
\end{align}
Then we multiply $\eqref{linearized problem of tilde-c}_1$ by $(-\Delta)\tilde{c}$ and integrate it over $\Omega\times(t_0,t)$ for any $t_0< t\le t_0+T$.
By Lemma \ref{G-N}, we have
\begin{align*}
&\frac{1}{2}\int_{\Omega}(|\nabla \tilde{c}(x,t)|^2-|\nabla \tilde{c}(x,t_0)|^2)\,\mathrm{d}x+\int_{t_0}^t\int_{\Omega}(|\Delta \tilde{c}|^2+(1-\eta+n_+)|\nabla \tilde{c}|^2)\,\mathrm{d}x\mathrm{d}s
   \nonumber
\\
& =   \int_{t_0}^t\int_{\Omega}n_+\nabla\cdot(\tilde{c}\nabla \tilde{c})\,\mathrm{d}x\mathrm{d}s+\int_{t_0}^t\int_{\Omega}u\cdot\nabla \tilde{c}\Delta\tilde{c}\,\mathrm{d}x\mathrm{d}s
-2\eta\int_{t_0}^t\int_{\Omega}\nabla g_1\cdot\nabla \tilde{c}\Delta\tilde{c}\,\mathrm{d}x\mathrm{d}s
\\
& \hspace{12pt} -\int_{t_0}^t\int_{\Omega}\eta(\eta|\nabla g_1|^2+\Delta g_1-u\nabla g_1-g_{1t})\tilde{c}\Delta\tilde{c}\,\mathrm{d}x\mathrm{d}s
\\
& \hspace{12pt} -\int_{t_0}^t\int_{\Omega}(\Delta g_2-u\nabla g_2-n_+g_2-g_{2t}-(1-\eta)g_2)e^{-\eta g_1}\Delta\tilde{c}\,\mathrm{d}x\mathrm{d}s
\\
& \le  \int_{t_0}^t\|n\|_{L^2}(\|\nabla \tilde{c}\|^2_{L^4}+\|\tilde{c}\|_{L^{\infty}}\|\Delta \tilde{c}\|_{L^2})\,\mathrm{d}s+\int_{t_0}^t\|\nabla u\|_{L^2}\|\nabla\tilde{c}\|^2_{L^4}\,\mathrm{d}x\mathrm{d}s
\\
&\hspace{12pt} +C\int_{t_0}^t\|\nabla \tilde{c}\|_{L^2}\|\Delta \tilde{c}\|_{L^2}\,\mathrm{d}s+\int_{t_0}^t(\|n\|_{L^2}+\|u\|_{L^2}+C)\|\Delta \tilde{c}\|_{L^2}\,\mathrm{d}s
\\
& \le  C\int_{t_0}^t(\|n\|_{L^2}+\|\nabla u\|_{L^2})(\|\tilde{c}\|_{L^{\infty}}\|\Delta \tilde{c}\|_{L^2}+\|\tilde{c}\|_{L^2})\,\mathrm{d}s
\\
&\hspace{12pt} +\frac{1}{4}\int_{t_0}^t\|\Delta \tilde{c}\|^2_{L^2}\,\mathrm{d}s+C\int_{t_0}^t(\|n\|_{L^2}^2+\|u\|^2_{L^2}+\|\nabla \tilde{c}\|^2_{L^2})\,\mathrm{d}s
\\
& \le  \frac{1}{2}\int_{t_0}^t\|\Delta \tilde{c}\|^2_{L^2}\,\mathrm{d}s+C\int_{t_0}^t(\|n\|_{L^2}^2+\|u\|^2_{H^1}+\|\tilde{c}\|^2_{H^1})\,\mathrm{d}s.
\end{align*}
Then, combining \eqref{priori estimate for fourth order system-2}--\eqref{priori estimate for fourth order system-3} with Lemma \ref{periodic gronwall inequality-1}, we obtain
\begin{align}\label{priori estimate for fourth order system-4}
\sup_{t}\|\nabla\tilde{c}(t)\|^2_{L^2}+\int_{0}^T\|\Delta \tilde{c}(s)\|^2_{L^2}\,\mathrm{d}s
\le  C\int_{0}^T\|n\|^2_{L^2}\,\mathrm{d}s.
\end{align}
Combining \eqref{priori estimate for fourth order system-1} with \eqref{priori estimate for fourth order system-4} and noting $s>2$, we have
\begin{align*}
&\sup_{t}\int_{\Omega}|n(x,t)|^2\,\mathrm{d}x+\int_{0}^T\int_{\Omega}(2\delta|\Delta n|^2+\varepsilon|n|^{s+2}+\mu|n|^3)\,\mathrm{d}x\mathrm{d}s
\\
&\le  C\int_{0}^T(\|n\|^4_{L^4}+\|n\|^2_{L^2})\,\mathrm{d}s+C
\\
&\le  \frac{\varepsilon}{2}\int_{0}^T\int_{\Omega}|n|^{s+2}\,\mathrm{d}x\mathrm{d}s+C,
\end{align*}
which implies
\begin{align}\label{priori estimate for fourth order system-5}
\sup_{t}\|n(t)\|^2_{L^2}&+\int_{t_0}^t\int_{\Omega}m(|n|+\varepsilon)^{m-1}|\nabla n|^2\,\mathrm{d}x\mathrm{d}s   \nonumber
\\
&+\int_{0}^T(\delta\|\Delta n\|_{L^2}^2+\varepsilon\|n\|_{L^{s+2}}^{s+2}+\mu\|n\|_{L^3}^3)\,\mathrm{d}s\le C,
\end{align}
where $C$ is independent of $\delta$ and $A$, but depends on $\varepsilon$.
Recalling \eqref{priori estimate for fourth order system-2} and Lemma \ref{solution for linearized problem of tilde-c}, we also have
\begin{align}\label{priori estimate for fourth order system-6}
&\sup_{t}(\|u(t)\|^2_{H^1}+\|\tilde{c}(t)\|^2_{H^2})+\int_{0}^T(\|\nabla u(s)\|^2_{H^1}+\|\nabla \Delta \tilde{c}(s)\|^2_{L^2})\,\mathrm{d}s\le C,
\end{align}
where $C$ is also independent of $\delta$ and $A$, but depends on $\varepsilon$.
Take $\hat{n}=n$ in \eqref{linearized problem of n-fourth order}.
Multiplying $\eqref{linearized problem of n-fourth order}_1$ by $(-\Delta)n$, integrating it over $\Omega\times(t_0,t)$ for
any $t_0< t\le t_0+T$, and using \eqref{priori estimate for fourth order system-5}--\eqref{priori estimate for fourth order system-6} and $s>2(m-1)$, we see that
\begin{align*}
&\frac{1}{2}\int_{\Omega}(|\nabla n(x,t)|^2-|\nabla n(x,t_0)|^2)\,\mathrm{d}x
+\int_{t_0}^t\int_{\Omega}(\delta|\nabla\Delta n|^2+\varepsilon(s+1)|n|^{s}|\nabla n|^2+2\mu|n||\nabla n|^2+A|\nabla n|^2)\,\mathrm{d}x\mathrm{d}s
\\
&  =  \int_{t_0}^t\int_{\Omega}m(|n|+\varepsilon)^{m-1}\nabla n\cdot\nabla\Delta n\,\mathrm{d}x\mathrm{d}s+\int_{t_0}^t\int_{\Omega}u\cdot\nabla n\Delta n\,\mathrm{d}x\mathrm{d}s
\\
& \hspace{12pt} -\eta\chi\int_{t_0}^t\int_{\Omega}(e^{g_1}n_+\nabla \tilde{c}+n_+\nabla g_2)\cdot\nabla\Delta n\,\mathrm{d}x\mathrm{d}s
+\eta\chi\int_{t_0}^t\int_{\Omega}\nabla\cdot(e^{g_1}n_+\tilde{c}\nabla g_1)\Delta n\,\mathrm{d}x\mathrm{d}s
\\
& \hspace{12pt} -\int_{t_0}^t\int_{\Omega}\eta(\mu a |n|+A|n|+g)\Delta n\,\mathrm{d}x\mathrm{d}s
\\
& \le  \frac{m^2}{\delta}\int_{t_0}^t\int_{\Omega}(|n|+\varepsilon)^{2(m-1)}|\nabla n|^2\,\mathrm{d}x\mathrm{d}s+\frac{\delta}{2}\int_{t_0}^t\|\nabla\Delta n\|_{L^2}^2\,\mathrm{d}s
+C\int_{t_0}^t\|u\|_{L^6}\|\nabla n\|_{L^3}\|\Delta n\|_{L^2}\,\mathrm{d}s
\\
& \hspace{12pt}  +C\int_{t_0}^t(\|n\nabla\tilde{c}\|_{L^2}^2+\|n\|_{L^2}^2+\|\Delta n\|_{L^2}^2)\,\mathrm{d}s+C\int_{t_0}^t(\|n\tilde{c}\|_{L^2}+\|\nabla(n\tilde{c})\|_{L^2})\|\Delta n\|_{L^2}\,\mathrm{d}s
\\
& \hspace{12pt}  +\eta A\int_{t_0}^t\|\nabla n\|_{L^2}^2\,\mathrm{d}s
\\
& \le  \frac{\varepsilon}{2}\int_{t_0}^t\int_{\Omega}|n|^{s}|\nabla n|^2\,\mathrm{d}x\mathrm{d}s+C\int_{t_0}^t\int_{\Omega}|\nabla n|^2\,\mathrm{d}x\mathrm{d}s
+\frac{\delta}{2}\int_{t_0}^t\|\nabla\Delta n\|_{L^2}^2\,\mathrm{d}s+\eta A\int_{t_0}^t\|\nabla n\|_{L^2}^2\,\mathrm{d}s
\\
&\hspace{12pt}  +C\int_{t_0}^t\|u\|_{H^1}\|\nabla n\|_{H^1}\|\Delta n\|_{L^2}\,\mathrm{d}s+C\int_{t_0}^t(\|n\nabla\tilde{c}\|_{L^2}^2+\|n\|_{L^2}^2+\|\Delta n\|_{L^2}^2)\,\mathrm{d}s.
\end{align*}
Since
$$
\int_{t_0}^t\|u\|_{H^1}\|\nabla n\|_{H^1}\|\Delta n\|_{L^2}\,\mathrm{d}s
\le \sup_{t}\|u(t)\|^2_{H^1}\int_{t_0}^t(\|\nabla n\|_{H^1}^2+\|\Delta n\|_{L^2}^2)\,\mathrm{d}s,
$$
$$
\int_{t_0}^t\|n\nabla\tilde{c}\|_{L^2}^2\,\mathrm{d}s
\le C \int_{t_0}^t\|n\|_{L^4}^2\|\nabla\tilde{c}\|_{L^4}^2\,\mathrm{d}s
\le C\sup_{t}\|\tilde{c}(t)\|^2_{H^2}\int_{t_0}^t\|n\|_{H^1}^2\,\mathrm{d}s,
$$
we have
\begin{align*}
&\int_{\Omega}(|\nabla n(x,t)|^2-|\nabla n(x,t_0)|^2)\,\mathrm{d}x
\\
& \hspace{12pt}  +\int_{t_0}^t\int_{\Omega}(\delta|\nabla\Delta n|^2+\varepsilon(s+1)|n|^{s}|\nabla n|^2+\mu|n||\nabla n|^2)\,\mathrm{d}x\mathrm{d}s
\le C\int_{t_0}^t\|n\|_{H^2}^2\,\mathrm{d}s.
\end{align*}
By Lemma \ref{periodic gronwall inequality-2} and \eqref{priori estimate for fourth order system-5}, we obtain
\begin{align}\label{priori estimate for fourth order system-7}
&\sup_{t}\|\nabla n(t)\|^2_{L^2}+\int_{t_0}^t\int_{\Omega}(\delta|\nabla\Delta n|^2+\varepsilon(s+1)|n|^{s}|\nabla n|^2+\mu|n||\nabla n|^2)\,\mathrm{d}x\mathrm{d}s
\le C,
\end{align}
where $C$ depends on $\varepsilon$ and $\delta$, and is independent of $A$.
By \eqref{priori estimate for fourth order system-5} and \eqref{priori estimate for fourth order system-7}, we
complete the proof of Lemma \ref{priori estimate for fourth order system}.
\hfill$\Box$

\vbox to 3mm{}

By Lemma \ref{solution for linearized problem of n-fourth order}, we see that $\mathcal{F}$ is a compact operator. Furthermore, by \eqref{priori estimate for fourth order system-1}, it
is easy to see that
$$
\mathcal{F}(n, 0)=0.
$$
Combining this equality with Lemma \ref{priori estimate for fourth order system}, and applying Leray-Schauder's fixed point theorem,
we see that the operator $\mathcal{F}(\cdot, 1)$ has a fixed point in $\mathcal{E}$, that is, the following problem admits a solution $(n, \tilde{c}, u)$,
\begin{equation}\label{fourth order problem-A}
 \left\{
\begin{array}{l}
n_t-m\nabla\cdot((|n|+\varepsilon)^{m-1}\nabla n)+\delta\Delta^2n+\varepsilon|n|^sn+u\cdot\nabla n+An
\\[1mm]
\hspace{35pt} =-\chi\nabla\cdot(e^{g_1}n_+\nabla \tilde{c}+n_+\tilde{c}\nabla e^{g_1}+n_+\nabla g_2)+(\mu a+A)|n|-\mu|n|n+g, \\[2mm]
\tilde{c}_t-\Delta \tilde{c}+(u-2\nabla g_1)\cdot\nabla \tilde{c}=(|\nabla g_1|^2+\Delta g_1-n_+-u\nabla g_1-g_{1t})\tilde{c}
\\[1mm]
\hspace{35pt} +(\Delta g_2-u\nabla g_2-n_+g_2-g_{2t})e^{-g_1}, \\[2mm]
u_t= \Delta u  -\nabla \pi +n\nabla \varphi, \\[2mm]
\nabla\cdot u =0,\\[2mm]
\frac{\partial n}{\partial \nu}|_{\partial\Omega}=\frac{\partial \Delta n}{\partial \nu}|_{\partial\Omega}
=\frac{\partial \tilde{c}}{\partial \nu}|_{\partial\Omega}=u|_{\partial\Omega}=0,
 \end{array} \right.
 \end{equation}
where $s>\max\{2(m-1), 2\}$, $\delta, \varepsilon>0$.
Taking advantage of Lemma \ref{solution for linearized problem of n-fourth order}, we have the following proposition.

\begin{proposition}\label{existence for fourth order system-A}
Assume $s>\max\{2(m-1),2\}$, and $a, g, \nabla\varphi\in L^{\infty}_T(Q)$. Then
the problem \eqref{fourth order problem-A} admits a strong time periodic solution $(n, \tilde{c}, u)$ with
\begin{align*}
&u\in  L^{\infty}_T(\mathbb{R}^+, H^1_{\sigma}(\Omega))\cap L^{2}_T(\mathbb{R}^+, H^2_{\sigma}(\Omega)), \quad u_t\in L^{2}_T(\mathbb{R}^+, L^2_{\sigma}(\Omega)),
\\
&\tilde{c}\in L^{\infty}_T(\mathbb{R}^+, H^2(\Omega))\cap L^{2}_T(\mathbb{R}^+, H^3(\Omega)), \quad \tilde{c}_t\in L^{2}_T(\mathbb{R}^+, L^2(\Omega)),
\\
&n\in L^{\infty}_T(\mathbb{R}^+, H^2(\Omega))\cap L^{2}_T(\mathbb{R}^+, H^4(\Omega)), \quad n_t\in L^{2}_T(\mathbb{R}^+, L^2(\Omega)).
\end{align*}
\end{proposition}

\vspace{5mm}

\setcounter{equation}{0}
\section{Weak time periodic solutions}

In this section, we  prove the existence of weak
time periodic solutions for the problem \eqref{homogeneous problem}.
As in \cite{Jin2020DCDS}, the proof is based on two level approximation schemes (corresponding to $\delta$ and $\varepsilon$, respectively).
We first consider the approximation corresponding to that $\delta\to0$.

\begin{lemma}\label{energy inequalities for fourth order system}
Assume $m>1$, $s>\max\{2(m-1),2\}$, and $a, g, \nabla\varphi\in L^{\infty}_T(Q)$. Let $(n, \tilde{c}, u)$ be the strong time periodic solution of \eqref{fourth order problem-A} obtained in Proposition \ref{existence for fourth order system-A}.
Then we have
\begin{align}\label{energy inequalities for fourth order system-1}
&\sup_{t}\int_{\Omega}(\delta|\nabla n|^2+(|n|+\varepsilon)^{m+1})\,\mathrm{d}x+\int_{t_0}^t\int_{\Omega}|\delta\nabla\Delta n-m(|n|+\varepsilon)^{m-1}\nabla n|^2\,\mathrm{d}x\mathrm{d}s   \nonumber
\\
& \hspace{35pt}  +\int_{t_0}^t\int_{\Omega}(|n|+\varepsilon)^{m}(\varepsilon|n|^{s+1}+\mu n^2)\,\mathrm{d}x\mathrm{d}s \le C,
\end{align}
where C is independent of $\delta$, and depends only on $\varepsilon, \Omega, T, a, \nabla\varphi$. Recalling $\eqref{fourth order problem-A}_1$, if
$s\le 5m-1$, we also have $n_t\in L^{\frac{6}{5}}(0,T;H^{-1}(\Omega))$.
\end{lemma}

\noindent{\bf Proof}.\
Multiplying $\eqref{fourth order problem-A}_1$ by $-\delta\Delta n+(|n|+\varepsilon)^{m}{\rm sgn}n$
and using \eqref{priori estimate for fourth order system-5} and \eqref{priori estimate for fourth order system-6}, we obtain
\begin{align*}
&\frac{\delta}{2}\int_{\Omega}(|\nabla n(x,t)|^2-|\nabla n(x,t_0)|^2)\,\mathrm{d}x
+\frac{1}{m+1}\int_{\Omega}((|n(x,t)|+\varepsilon)^{m+1}-(|n(x,t_0)|+\varepsilon)^{m+1})\,\mathrm{d}x
\\
& \hspace{12pt} +\int_{t_0}^t\int_{\Omega}|\delta\nabla\Delta n-m(|n|+\varepsilon)^{m-1}\nabla n|^2\,\mathrm{d}x\mathrm{d}s+\int_{t_0}^t\int_{\Omega}(\varepsilon\delta(s+1)|n|^s+2\mu|n|+A\delta)|\nabla n|^2\,\mathrm{d}x\mathrm{d}s
\\
& \hspace{12pt} +\int_{t_0}^t\int_{\Omega}(\varepsilon|n|^{s+1}+\mu n^2+A|n|)(|n|+\varepsilon)^{m}\,\mathrm{d}x\mathrm{d}s
\\
&  =  \delta\int_{t_0}^t\int_{\Omega}u\cdot\nabla n\Delta n\,\mathrm{d}x\mathrm{d}s-\chi\int_{t_0}^t\int_{\Omega}(e^{g_1}n_+\nabla \tilde{c}+n_+\nabla g_2)(\delta\nabla\Delta n-m(|n|+\varepsilon)^{m-1}\nabla n)\,\mathrm{d}x\mathrm{d}s
\\
& \hspace{12pt}  +\chi\int_{t_0}^t\int_{\Omega}\nabla\cdot(e^{g_1}n_+\tilde{c}\nabla g_1)(\delta\Delta n-(|n|+\varepsilon)^{m}{\rm sgn}n)\,\mathrm{d}x\mathrm{d}s
\\
& \hspace{12pt} -\int_{t_0}^t\int_{\Omega}(\mu a |n|+A|n|+g)(\delta\Delta n-(|n|+\varepsilon)^{m}{\rm sgn}n)\,\mathrm{d}x\mathrm{d}s
\\
& \le  \delta\int_{t_0}^t\|u\|_{L^6}\|\nabla n\|_{L^3}\|\Delta n\|_{L^2}\,\mathrm{d}s
+\frac{1}{2}\int_{t_0}^t\int_{\Omega}|\delta\nabla\Delta n-m(|n|+\varepsilon)^{m-1}\nabla n|^2\,\mathrm{d}x\mathrm{d}s
\\
& \hspace{12pt}  +C\int_{t_0}^t(\|n\|^2_{L^3}\|\nabla\tilde{c}\|_{L^6}^2+\|n\|_{H^1}^2)\,\mathrm{d}s+2\delta\int_{t_0}^t\|\Delta n\|_{L^2}^2\,\mathrm{d}s
+C\int_{t_0}^t\int_{\Omega}(|n|+\varepsilon)^{2m}\,\mathrm{d}x\mathrm{d}s
\\
& \hspace{12pt}  +\delta A\int_{t_0}^t\|\nabla n\|_{L^2}^2\,\mathrm{d}s+A\int_{t_0}^t\int_{\Omega}|n|(|n|+\varepsilon)^{m}\,\mathrm{d}x\mathrm{d}s
+\int_{t_0}^t\int_{\Omega}(|n|+\varepsilon)^{m}\,\mathrm{d}x\mathrm{d}s
\\
& \le  \delta\int_{t_0}^t\|\nabla u\|_{L^2}\|\nabla n\|^{\frac{1}{4}}_{L^2}\|\Delta n\|^{\frac{7}{4}}_{L^2}\,\mathrm{d}s
+\frac{1}{2}\int_{t_0}^t\int_{\Omega}|\delta\nabla\Delta n-m(|n|+\varepsilon)^{m-1}\nabla n|^2\,\mathrm{d}x\mathrm{d}s
\\
&\hspace{12pt}  +C\int_{t_0}^t(\|n\|_{L^3}^3+\|\nabla\tilde{c}\|_{H^1}^6+\|n\|_{H^1}^2)\,\mathrm{d}s+2\delta\int_{t_0}^t\|\Delta n\|_{L^2}^2\,\mathrm{d}s
+C\int_{t_0}^t\int_{\Omega}(|n|+\varepsilon)^{2m}\,\mathrm{d}x\mathrm{d}s
\\
& \hspace{12pt}  +\delta A\int_{t_0}^t\|\nabla n\|_{L^2}^2\,\mathrm{d}s+A\int_{t_0}^t\int_{\Omega}|n|(|n|+\varepsilon)^{m}\,\mathrm{d}x\mathrm{d}s
+\int_{t_0}^t\int_{\Omega}(|n|+\varepsilon)^{m}\,\mathrm{d}x\mathrm{d}s
\\
& \le  3\delta\int_{t_0}^t\|\Delta n\|_{L^2}^2\,\mathrm{d}s
+\frac{1}{2}\int_{t_0}^t\int_{\Omega}|\delta\nabla\Delta n-m(|n|+\varepsilon)^{m-1}\nabla n|^2\,\mathrm{d}x\mathrm{d}s
\\
& \hspace{12pt}  +C\int_{t_0}^t(\|\nabla\tilde{c}\|_{H^1}^2+\|n\|_{L^3}^3+\|n\|_{H^1}^2+\|u\|_{H^1}^2)\,\mathrm{d}s
+\frac{\varepsilon}{2}\int_{t_0}^t\int_{\Omega}|n|^{m+s+1}\,\mathrm{d}x\mathrm{d}s
\\
&\hspace{12pt}  +\delta A\int_{t_0}^t\|\nabla n\|_{L^2}^2\,\mathrm{d}s+A\int_{t_0}^t\int_{\Omega}|n|(|n|+\varepsilon)^{m}\,\mathrm{d}x\mathrm{d}s+C.
\end{align*}
Then we have
\begin{align*}
&\frac{\delta}{2}\int_{\Omega}(|\nabla n(x,t)|^2-|\nabla n(x,t_0)|^2)\,\mathrm{d}x
+\frac{1}{m+1}\int_{\Omega}((|n(x,t)|+\varepsilon)^{m+1}-(|n(x,t_0)|+\varepsilon)^{m+1})\,\mathrm{d}x
\\
& \hspace{12pt} +\frac{1}{2}\int_{t_0}^t\int_{\Omega}|\delta\nabla\Delta n-m(|n|+\varepsilon)^{m-1}\nabla n|^2\,\mathrm{d}x\mathrm{d}s+\int_{t_0}^t\int_{\Omega}(\varepsilon\delta(s+1)|n|^s+2\mu|n|)|\nabla n|^2\,\mathrm{d}x\mathrm{d}s
\\
& \hspace{12pt} +\frac{1}{2}\int_{t_0}^t\int_{\Omega}(\varepsilon|n|^{s+1}+\mu n^2)(|n|+\varepsilon)^{m}\,\mathrm{d}x\mathrm{d}s
\\
& \le  3\delta\int_{t_0}^t\|\Delta n\|_{L^2}^2\,\mathrm{d}s+C\int_{t_0}^t(\|\nabla\tilde{c}\|_{H^1}^2+\|n\|_{L^3}^3+\|n\|_{H^1}^2+\|u\|_{H^1}^2)\,\mathrm{d}s+C.
\end{align*}
Combining the above inequality with \eqref{priori estimate for fourth order system-5}--\eqref{priori estimate for fourth order system-6} and
applying Lemma \ref{periodic gronwall inequality-2}, we obtain \eqref{energy inequalities for fourth order system-1}. By \eqref{energy inequalities for fourth order system-1},
it is easy to see that $\frac{6(s+1)}{5}\le s+1+m$ and then $\varepsilon|n|^{s+1}\in L^{\frac{6}{5}}(0,T; L^{\frac{6}{5}}(\Omega))$ if $s\le5m-1$.
Noting that $L^{\frac{6}{5}}(\Omega)\hookrightarrow H^{-1}(\Omega)$, we obtain $n_t\in L^{\frac{6}{5}}(0,T; H^{-1}(\Omega))$.
\hfill$\Box$

\vbox to 3mm{}

Let $(u_{\varepsilon\delta}, \tilde{c}_{\varepsilon\delta}, n_{\varepsilon\delta})$
be a time periodic solution of the problem \eqref{fourth order problem-A} satisfying \eqref{priori estimate for fourth order system-5},
\eqref{priori estimate for fourth order system-6} and \eqref{energy inequalities for fourth order system-1}.
Taking $\delta\to0$,
we obtain that (if necessary, we may choose a subsequence)
\begin{equation*}
\begin{array}{l@{\hspace{-1pt}}l@{\hspace{15pt}}l}
&u_{\varepsilon\delta}\to u_{\varepsilon}, & {\rm in}\ \ L^6(Q_T),
\\[2mm]
&\tilde{c}_{\varepsilon\delta}\to \tilde{c}_{\varepsilon}, & {\rm in}\ \ C(\overline{Q_T}),
\\[2mm]
&u_{\varepsilon\delta}\rightharpoonup u_{\varepsilon},  \ \tilde{c}_{\varepsilon\delta}\rightharpoonup \tilde{c}_{\varepsilon}, & {\rm in}\ \  W_2^{2,1}(Q_T),
\\[2mm]
&n_{\varepsilon\delta}\to n_{\varepsilon}, & {\rm in}\ \ L^p(Q_T)\ \ {\rm for\ any}\ \ p<s+m+1,
\\[2mm]
&\nabla n_{\varepsilon\delta}\rightharpoonup \nabla n_{\varepsilon},  \ \delta\Delta n_{\varepsilon\delta}\rightarrow 0, & {\rm in}\ \ L^2(Q_T),
\\[2mm]
&n_{\varepsilon\delta}\rightharpoonup n_{\varepsilon}, \ \tilde{c}_{\varepsilon\delta}\to \tilde{c}_{\varepsilon}, & {\rm in}\ \ L^2(\partial\Omega\times(0, T)).
\end{array}
\end{equation*}
Then $(u_{\varepsilon}, \tilde{c}_{\varepsilon}, n_{\varepsilon})$
is the solution of the following problem
\begin{equation}\label{second order problem-A}
 \left\{
\begin{array}{l}
n_t-m\nabla\cdot((|n|+\varepsilon)^{m-1}\nabla n)+\varepsilon|n|^sn+u\cdot\nabla n+An
\\[1mm]
\hspace{35pt} =-\chi\nabla\cdot(e^{g_1}n_+\nabla \tilde{c}+n_+\tilde{c}\nabla e^{g_1}+n_+\nabla g_2)+(\mu a+A)|n|-\mu|n|n+g, \\[2mm]
\tilde{c}_t-\Delta \tilde{c}+(u-2\nabla g_1)\cdot\nabla \tilde{c}=(|\nabla g_1|^2+\Delta g_1-n_+-u\nabla g_1-g_{1t})\tilde{c}
\\[1mm]
\hspace{35pt} +(\Delta g_2-u\nabla g_2-n_+g_2-g_{2t})e^{-g_1}, \\[2mm]
u_t= \Delta u  -\nabla \pi +n\nabla \varphi, \\[2mm]
\nabla\cdot u =0,\\[2mm]
\frac{\partial n}{\partial \nu}|_{\partial\Omega}
=\frac{\partial \tilde{c}}{\partial \nu}|_{\partial\Omega}=u|_{\partial\Omega}=0,
 \end{array} \right.
 \end{equation}
such that $(u_{\varepsilon}, \tilde{c}_{\varepsilon}, n_{\varepsilon})$ satisfies \eqref{priori estimate for fourth order system-5} and \eqref{priori estimate for fourth order system-6}.

Next, we consider the second level approximation. For this purpose, we need the following energy estimates
for the problem \eqref{second order problem-A} independent of $\varepsilon$.

\begin{proposition}\label{estimate for second order system}
Assume $m>1$, $\max\{2(m-1),2\}<s\le5m-1$, $g\ge0$ and $a, g, \nabla\varphi\in L^{\infty}_T(Q)$.
If $A>0$ is sufficiently large, the problem \eqref{second order problem-A} admits a time periodic solution $(n_{\varepsilon}, \tilde{c}_{\varepsilon}, u_{\varepsilon})$
such that $n_{\varepsilon}, c_{\varepsilon}\ge0$ and
\begin{align*}
&\sup_{t}(\|u_{\varepsilon}(t)\|_{W^{1,\infty}}+\|n_{\varepsilon}(t)\|_{L^\infty}+\|\tilde{c}_{\varepsilon}(t)\|_{W^{1,\infty}})\le C,
\\
&\sup_{t}(\|u_{\varepsilon}(t)\|^2_{H^1}+\|\tilde{c}_{\varepsilon}(t)\|^2_{H^1})+\int_{0}^T(\|u_{\varepsilon}\|^2_{H^2}+\|\tilde{c}_{\varepsilon}\|^2_{H^2}+\|u_{\varepsilon t}\|^2_{L^2}+\|\tilde{c}_{\varepsilon t}\|^2_{L^2})\,\mathrm{d}s\le C,
\\
&\sup_{t}(\|n_{\varepsilon}\|_{L^1}+\|n_{\varepsilon}\ln n_{\varepsilon}\|_{L^1})+\int_{0}^T\int_{\Omega}(n_{\varepsilon}+\varepsilon)^{m-1}|\nabla \sqrt{n_{\varepsilon}}|^2\,\mathrm{d}x\mathrm{d}s
\le C,
\\
&\int_{0}^T\left(\|u_t\|^p_{L^p}+\|u\|^p_{W^{2,p}}+\|\tilde{c}_t\|^p_{L^p}+\|\tilde{c}\|^p_{W^{2,p}}\right)\,\mathrm{d}s\le C \ \ {\rm for\ any\ } p>1,
\\
&\sup_t\int_{\Omega}|\nabla (n_{\varepsilon}+\varepsilon)^m|^2\,\mathrm{d}x
+\int_{0}^T\int_{\Omega}(n_{\varepsilon}+\varepsilon)^{m-1}\left|\frac{\partial n_{\varepsilon}}{\partial t}\right|^2\,\mathrm{d}x\mathrm{d}s \le C,
\end{align*}
where $C$ is independent of $\varepsilon$.
\end{proposition}

The proof of Proposition \ref{estimate for second order system} is give by the following lemmas.
For simplicity, we will denote the limit $(u_{\varepsilon}, \tilde{c}_{\varepsilon}, n_{\varepsilon})$ by $(u, \tilde{c}, n)$.

\begin{lemma}\label{L1 estimate for second order system}
Assume $m>1$, $\max\{2(m-1),2\}<s\le5m-1$, $g\ge0$ and $a, g, \nabla\varphi\in L^{\infty}_T(Q)$. Let $(n, \tilde{c}, u)$ be a periodic solution of \eqref{second order problem-A}.
Then for sufficiently large $A$, we have $n, c\ge0$ and
\begin{align}
&\sup_{t}\|\tilde{c}(t)\|_{L^\infty}\le C,\label{L-infty estimate of u and c}
\\
&\sup_{t}(\|u(t)\|^2_{H^1}+\|\tilde{c}(t)\|^2_{H^1})+\int_{0}^T(\|u\|^2_{H^2}+\|\tilde{c}\|^2_{H^2}+\|u_t\|^2_{L^2}+\|\tilde{c}_t\|^2_{L^2})\,\mathrm{d}s\le C, \label{L1 estimate for second order system-1}
\\\label{L1 estimate for second order system-2}
&\sup_{t}(\|n\|_{L^1}+\|n\ln n\|_{L^1})+\int_{0}^T\int_{\Omega}(n+\varepsilon)^{m-1}|\nabla \sqrt{n}|^2\,\mathrm{d}x\mathrm{d}s+\mu\int_{0}^T\|n\|^2_{L^2}\,\mathrm{d}s
\le C,
\end{align}
where $C$ is independent of $\varepsilon$.
\end{lemma}

\noindent{\bf Proof}.\
By \cite[Lemma 3.2]{Jin2020PRSE}, we have
$$
0\le c\le\|g_2\|_{L^{\infty}}.
$$
Next, we show $n\ge0$ by
examining the set $J(t)=\{x\in\Omega; n(x,t)<0\}$. Assume that $J(t)$ is a differentiable submanifold and $\frac{\partial n}{\partial \nu'}$
denote the outward normal derivative of $n$ on $J(t)$.
It's easy to see that $n=0$ and $\frac{\partial n}{\partial \nu'}\ge0$ on $\partial\{J(t)\}\setminus\partial\Omega$,
and $n_+=0$ and $\frac{\partial n}{\partial \nu'}=0$ on $\partial\{J(t)\}\cap\partial\Omega$.
A direct integration of $\eqref{second order problem-A}_1$ on $J(t)\times(0,T)$ gives
\begin{align*}
0\ge&-m\int_0^T\int_{\partial\{J(t)\}}(|n|+\varepsilon)^{m-1}\frac{\partial n}{\partial \nu'}\,\mathrm{d}\Gamma\mathrm{d}s
+\int_0^T\int_{J(t)}(\varepsilon|n|^s+A)n\,\mathrm{d}x\mathrm{d}s
\\
&=\int_0^T\int_{J(t)}(\mu a+A)|n|\,\mathrm{d}x\mathrm{d}s+\int_0^T\int_{J(t)}(g-\mu|n|n)\,\mathrm{d}x\mathrm{d}s\ge0,
\end{align*}
where $A>0$ is sufficiently large.
It implies that
$$
\int_0^T\int_{J(t)}n\,\mathrm{d}x\mathrm{d}s=0,
$$
namely, $n\ge0$. While if $J(t)$ is not a regular submanifold, we can construct a sufficiently
smooth approximating sequence $(n_k, \tilde{c}_k, u_k)$ of $(n, \tilde{c}, u)$ such that the
corresponding approximating solutions $n_k$ satisfying that $n_k(\cdot,t)$ are continuously
differentiable. Thus, the sets $J_k(t)$ are measurable and $\partial J_k(t)$ are differentiable
submanifolds. Then the above result can be obtained by letting $k\to0$.

Similarly, in what follows, we still assume that the solution $n$ is sufficiently
smooth. Otherwise, the following estimates can be obtained by an approximating
process. Take $n_+=|n|=n$ in the problem \eqref{second order problem-A}.
Combining \eqref{second order problem-A} with \eqref{priori estimate for fourth order system-2}--\eqref{priori estimate for fourth order system-4}, we have
\begin{align}\label{L1 estimate for second order system-3}
\sup_{t}(\|u(t)\|^2_{H^1}+\|\tilde{c}(t)\|^2_{H^1})+\int_{0}^T&(\|u\|^2_{H^2}+\|\tilde{c}\|^2_{H^2}+\|u_t\|^2_{L^2}+\|\tilde{c}_t\|^2_{L^2})\,\mathrm{d}s  \nonumber
\\
&\le C\int_{0}^T\|n\|^2_{L^2}\,\mathrm{d}s+C,
\end{align}
where C is independent of $\varepsilon$.
In order to estimate the term $\int_{0}^T\|n\|^2_{L^2}\,\mathrm{d}s$, we integrate the equation $\eqref{second order problem-A}_1$
over $\Omega\times(t_0,t_0+t)$ with $t_0<t\le t_0+T$.
It is easy to obtain that
\begin{align*}
&\int_{\Omega}(n(x,t)-n(x,t_0))\,\mathrm{d}x+\int_{t_0}^t\int_{\Omega}(\mu n^2+\varepsilon n^{s+1})\,\mathrm{d}x\mathrm{d}s
\\
&= -\chi\int_{t_0}^t\int_{\partial\Omega}e^{g_1}n\tilde{c}\frac{\partial g_1}{\partial \nu}\,\mathrm{d}\Gamma\mathrm{d}s
+\int_{t_0}^t\int_{\Omega}(\mu an+g)\,\mathrm{d}x\mathrm{d}s
\\
& \le  C\int_{t_0}^t\int_{\partial\Omega}n\,\mathrm{d}\Gamma\mathrm{d}s
+C\int_{t_0}^t\int_{\Omega}n\,\mathrm{d}x\mathrm{d}s+C.
\end{align*}
Combining with Lemma \ref{periodic gronwall inequality-2}, we have
\begin{align}\label{L1 estimate for second order system-4}
\sup_{t}\|n\|_{L^1}+\mu\int_0^T\|n\|^2_{L^2}\,\mathrm{d}s
\le C\int_{0}^T\int_{\partial\Omega}n\,\mathrm{d}\Gamma\mathrm{d}s+C.
\end{align}
Obviously, the estimate of the term $\int_{0}^T\int_{\partial\Omega}n\,\mathrm{d}\Gamma\mathrm{d}s$  is a key point
to prove \eqref{L1 estimate for second order system-1} and \eqref{L1 estimate for second order system-2}.

Multiplying $\eqref{second order problem-A}_1$ by $1+\ln n$, and integrating it over $\Omega\times(t_0,t)$ with $t_0<t \le t_0+T$, we have
\begin{align}\label{L1 estimate for second order system-5}
&\int_{\Omega}(n(x,t)\ln n(x,t)-n(x,t_0)\ln n(x,t_0))\,\mathrm{d}x+4m\int_{t_0}^t\int_{\Omega}(n+\varepsilon)^{m-1}|\nabla \sqrt{n}|^2\,\mathrm{d}x\mathrm{d}s    \nonumber
\\
& \le \chi\int_{t_0}^t\int_{\Omega}(e^{g_1}\nabla \tilde{c}+e^{g_1}\tilde{c}\nabla g_1+\nabla g_2)\nabla n\,\mathrm{d}x\mathrm{d}s
-\chi\int_{t_0}^t\int_{\partial\Omega}e^{g_1}\tilde{c}\frac{\partial g_1}{\partial \nu}n(1+\ln n)\,\mathrm{d}\Gamma\mathrm{d}s    \nonumber
\\
& \hspace{12pt} +\int_{t_0}^t\int_{\Omega}(\mu an+g)(1+\ln n)\,\mathrm{d}x\mathrm{d}s    \nonumber
\\
& \le  C\int_{t_0}^t(\|\nabla \tilde{c}\|^2_{L^2}+\|\nabla n\|^2_{L^2}+\| n\|^2_{L^2})\,\mathrm{d}s
+C\int_{t_0}^t\int_{\partial\Omega}n(1+\ln n)_+\,\mathrm{d}\Gamma\mathrm{d}s+C,
\end{align}
where we use the fact that
\begin{align}\label{inequality of log}
1+\ln n\le\frac{1}{\sigma}e^{(\sigma-1)_+}n^\sigma
\end{align}
for any $n>0$ and $\sigma>0$.
Applying Lemma \ref{periodic gronwall inequality-1} and combining with \eqref{L1 estimate for second order system-3} and \eqref{L1 estimate for second order system-4},
we obtain
\begin{align}\label{L1 estimate for second order system-6}
&\sup_{t}\int_{\Omega}n(x,t)\ln n(x,t)\,\mathrm{d}x+4m\int_{0}^T\int_{\Omega}(n+\varepsilon)^{m-1}|\nabla \sqrt{n}|^2\,\mathrm{d}x\mathrm{d}s     \nonumber
\\
& \le  \frac{C}{\mu}\int_{0}^T\int_{\partial\Omega}n\,\mathrm{d}\Gamma\mathrm{d}s
+C\int_{0}^T\int_{\partial\Omega}n(1+\ln n)_+\,\mathrm{d}\Gamma\mathrm{d}s+C.
\end{align}
Taking $\sigma\in(0,1)$, adding \eqref{L1 estimate for second order system-4} and \eqref{L1 estimate for second order system-6} together
and recalling the embedding relation $W^{1,1}(\Omega)\hookrightarrow L^1(\partial\Omega)$ by Lemma \ref{Trace Imbedding}, yield
\begin{align*}
&\sup_{t}(\|n\|_{L^1}+\|n\ln n\|_{L^1})+\mu\int_0^T\|n\|^2_{L^2}\,\mathrm{d}s+4m\int_{0}^T\int_{\Omega}(n+\varepsilon)^{m-1}|\nabla \sqrt{n}|^2\,\mathrm{d}x\mathrm{d}s     \nonumber
\\
& \le  C\frac{1+\mu}{\mu}\int_{0}^T\int_{\partial\Omega}n\,\mathrm{d}\Gamma\mathrm{d}s
+C\int_{0}^T\int_{\partial\Omega}n(1+\ln n)_+\,\mathrm{d}\Gamma\mathrm{d}s+C  \nonumber
\\
&  =  C\frac{1+\mu}{\mu}\int_{t_0}^t\int_{\partial\Omega}n(2+\ln n-(1+\ln n))\,\mathrm{d}\Gamma\mathrm{d}s
+C\int_{0}^T\int_{\partial\Omega}n(1+\ln n)_+\,\mathrm{d}\Gamma\mathrm{d}s+C  \nonumber
\\
& \le C\frac{1+\mu}{\mu}\int_{t_0}^t\int_{\partial\Omega}n(2+\ln n)_+\,\mathrm{d}\Gamma\mathrm{d}s
-C\frac{1+\mu}{\mu}\int_{t_0}^t\int_{\partial\Omega}n(1+\ln n)_+\,\mathrm{d}\Gamma\mathrm{d}s  \nonumber
\\
& \hspace{12pt} +C\int_{0}^T\int_{\partial\Omega}n(1+\ln n)_+\,\mathrm{d}\Gamma\mathrm{d}s+C  \nonumber
\\
& \le  C\frac{1+\mu}{\mu}\int_{t_0}^t\int_{\Omega}n(2+\ln n)_+\,\mathrm{d}x\mathrm{d}s+C\frac{1+\mu}{\mu}\int_{t_0}^t\int_{\Omega}|\nabla (n(2+\ln n)_+)|\,\mathrm{d}x\mathrm{d}s+C  \nonumber
\\
& \le C\int_{t_0}^t\int_{\Omega}n^{1+\sigma}\,\mathrm{d}x\mathrm{d}s+C\int_{t_0}^t\int_{\Omega}|\nabla n|(2+\ln n)_+\,\mathrm{d}x\mathrm{d}s    \nonumber
\\
& \hspace{12pt} +C\int_{t_0}^t\int_{\Omega\cap {\rm supp}\,(2+\ln n)_+}|\nabla n|\,\mathrm{d}x\mathrm{d}s+C    \nonumber
\\
& \le C\int_{t_0}^t\int_{\Omega}n^{1+\sigma}\,\mathrm{d}x\mathrm{d}s+C\int_{t_0}^t\int_{\Omega}|\nabla n|(2+\ln n)_+\,\mathrm{d}x\mathrm{d}s    \nonumber
\\
&\hspace{12pt} +C\int_{t_0}^t\int_{\Omega\cap {\rm supp}\,(2+\ln n)_+}|\nabla n||1+\ln n|\,\mathrm{d}x\mathrm{d}s+C    \nonumber
\\
& \le C\int_{t_0}^t\int_{\Omega}n^{1+\sigma}\,\mathrm{d}x\mathrm{d}s+C\int_{t_0}^t\int_{\Omega}|\nabla n|n^{\frac{m-1+\sigma}{2}}\,\mathrm{d}x\mathrm{d}s+C  \nonumber
\\
& \le  \frac{\mu}{2}\int_0^T\|n\|^2_{L^2}\,\mathrm{d}s+m\int_{t_0}^t\int_{\Omega}|\nabla n|^2n^{m-2}\,\mathrm{d}x\mathrm{d}s+C,
\end{align*}
which implies that
\begin{align}\label{L1 estimate for second order system-7}
&\sup_{t}(\|n\|_{L^1}+\|n\ln n\|_{L^1})+\mu\int_0^T\|n\|^2_{L^2}\,\mathrm{d}s+\int_{0}^T\int_{\Omega}(n+\varepsilon)^{m-1}|\nabla \sqrt{n}|^2\,\mathrm{d}x\mathrm{d}s
\le C,
\end{align}
where $C$ depends on $\mu$ and is independent of $\varepsilon$.
Combining \eqref{L1 estimate for second order system-3} with \eqref{L1 estimate for second order system-7}, we obtain
\eqref{L1 estimate for second order system-1} and \eqref{L1 estimate for second order system-2}.
The proof of Lemma \ref{L1 estimate for second order system} is completed.
\hfill$\Box$

\vbox to 3mm{}

Since that $n\ge0$, the problem \eqref{second order problem-A} can be reduced to the following form.
\begin{equation}\label{second order problem-B}
 \left\{
\begin{array}{l}
n_t-\Delta(n+\varepsilon)^m+\varepsilon n^{s+1}+u\cdot\nabla n
\\[1mm]
\hspace{35pt} =-\chi\nabla\cdot(e^{g_1}n\nabla \tilde{c}+n\tilde{c}\nabla e^{g_1}+n\nabla g_2)+\mu n(a-n)+g, \\[2mm]
\tilde{c}_t-\Delta \tilde{c}+(u-2\nabla g_1)\cdot\nabla \tilde{c}=(|\nabla g_1|^2+\Delta g_1-n-u\nabla g_1-g_{1t})\tilde{c}
\\[1mm]
\hspace{35pt} +(\Delta g_2-u\nabla g_2-ng_2-g_{2t})e^{-g_1}, \\[2mm]
u_t= \Delta u  -\nabla \pi +n\nabla \varphi, \\[2mm]
\nabla\cdot u =0,\\[2mm]
\frac{\partial n}{\partial \nu}|_{\partial\Omega}
=\frac{\partial \tilde{c}}{\partial \nu}|_{\partial\Omega}=u|_{\partial\Omega}=0.
 \end{array} \right.
 \end{equation}

\begin{lemma}\label{2m base estimate}
Assume $m>1$, $\max\{2(m-1),2\}<s\le5m-1$, $g\ge0$ and $a, g, \nabla\varphi\in L^{\infty}_T(Q)$. Let $(n, \tilde{c}, u)$ be the periodic solution of \eqref{second order problem-B}.
Then for any $q\ge2(m-1)$, we have
\begin{align}\label{2m base estimate-1}
&\sup_t\int_{\Omega}n^{1+q}\,\mathrm{d}x+\int_{0}^T\int_{\Omega}\left(n^{m+q-2}|\nabla n|^2+n^{2+q}\right)\mathrm{d}x\mathrm{d}s
\le  C(q),
\end{align}
where $C(q)$ is independent of $\varepsilon$, but depends on $q$.
\end{lemma}

\noindent{\bf Proof}.\
Multiplying $\eqref{second order problem-B}_1$ by $n^{q}$ with $q\ge2(m-1)$, and integrating it over $\Omega\times(t_0,t)$ with $t_0<t \le t_0+T$, we have
\begin{align}\label{2m base estimate-2}
&\frac{1}{1+q}\int_{\Omega}(n^{1+q}(x,t)-n^{1+q}(x,t_0))\,\mathrm{d}x+\int_{t_0}^t\int_{\Omega}\left(qmn^{m+q-2}|\nabla n|^2+\mu n^{2+q}\right)\,\mathrm{d}x\mathrm{d}s    \nonumber
\\
& \le q\chi\int_{t_0}^t\int_{\Omega}n^{q}(e^{g_1}\nabla \tilde{c}+e^{g_1}\tilde{c}\nabla g_1+\nabla g_2)\cdot\nabla n\,\mathrm{d}x\mathrm{d}s
-\chi\int_{t_0}^t\int_{\partial\Omega}e^{g_1}\tilde{c}\frac{\partial g_1}{\partial \nu}n^{1+q}\,\mathrm{d}\Gamma\mathrm{d}s    \nonumber
\\
& \hspace{12pt} +\int_{t_0}^t\int_{\Omega}(\mu an+g)n^{q}\,\mathrm{d}x\mathrm{d}s    \nonumber
\\
&  \le   \frac{q}{4}\int_{t_0}^t\int_{\Omega}n^{q+m-2}|\nabla n|^2\,\mathrm{d}x\mathrm{d}s
+C\int_{t_0}^t\int_{\Omega}n^{q-m+2}|\nabla \tilde{c}|^2\,\mathrm{d}x\mathrm{d}s    \nonumber
\\
& \hspace{12pt}    +C\int_{t_0}^t\int_{\Omega}(n^{q}+n^{q+1}+n^{q-m+2})\,\mathrm{d}x\mathrm{d}s+C\int_{t_0}^t\int_{\partial\Omega}n^{1+q}\,\mathrm{d}\Gamma\mathrm{d}s    \nonumber
\\
& \le  \frac{q}{4}\int_{t_0}^t\int_{\Omega}n^{q+m-2}|\nabla n|^2\,\mathrm{d}x\mathrm{d}s+\frac{\mu}{4}\int_{t_0}^t\int_{\Omega}n^{q+2}\,\mathrm{d}x\mathrm{d}s
+C\int_{t_0}^t\int_{\Omega}|\nabla\tilde{c}|^{\frac{2(q+2)}{m}}\,\mathrm{d}x\mathrm{d}s    \nonumber
\\
& \hspace{12pt} +C\int_{t_0}^t\int_{\partial\Omega}n^{1+q}\,\mathrm{d}\Gamma\mathrm{d}s+C.
\end{align}
Recalling the embedding relation $W^{1,1}(\Omega)\hookrightarrow L^1(\partial\Omega)$ by Lemma \ref{Trace Imbedding} yields
\begin{align}\label{2m base estimate-3}
&C\int_{t_0}^t\int_{\partial\Omega}n^{1+q}\,\mathrm{d}\Gamma\mathrm{d}s    \nonumber
\\
& \le C\int_{t_0}^t\int_{\Omega}n^{1+q}\,\mathrm{d}x\mathrm{d}s+Cq\int_{t_0}^t\int_{\Omega}n^{q}|\nabla n|\,\mathrm{d}x\mathrm{d}s    \nonumber
\\
& \le \frac{\mu}{4}\int_{t_0}^t\int_{\Omega}n^{2+q}\,\mathrm{d}x\mathrm{d}s+C\int_{t_0}^t\int_{\Omega}n^{q}\,\mathrm{d}x\mathrm{d}s    \nonumber
\\
& \hspace{12pt} +\frac{q}{4}\int_{t_0}^t\int_{\Omega}n^{m+q-2}|\nabla n|^2\,\mathrm{d}x\mathrm{d}s
+Cq\int_{t_0}^t\int_{\Omega}n^{q-m+2}\,\mathrm{d}x\mathrm{d}s.
\end{align}
Substituting \eqref{2m base estimate-3} into \eqref{2m base estimate-2} and combining with Lemma \ref{periodic gronwall inequality-1} give
\begin{align}\label{2m base estimate-4}
&\sup_t\int_{\Omega}n^{1+q}(x,t)\,\mathrm{d}x+\int_{0}^T\int_{\Omega}\left(qmn^{m+q-2}|\nabla n|^2+\mu n^{2+q}\right)\,\mathrm{d}x\mathrm{d}s    \nonumber
\\
& \le  C\int_{0}^T\int_{\Omega}|\nabla\tilde{c}|^{\frac{2(q+2)}{m}}\mathrm{d}x\mathrm{d}s+C    \nonumber
\\
& \le  C\sup_t\|\tilde{c}\|^{\frac{q+2}{m}}_{L^{\infty}}\int_{0}^T\|\Delta\tilde{c}\|^{\frac{q+2}{m}}_{L^{\frac{q+2}{m}}}\,\mathrm{d}s+C\sup_t\|\tilde{c}\|^{\frac{2(q+2)}{m}}_{L^{\infty}}+C   \nonumber
\\
& \le  C\int_{0}^T\|\Delta\tilde{c}\|^{\frac{q+2}{m}}_{L^{\frac{q+2}{m}}}\,\mathrm{d}s+C,
\end{align}
where $C$ is independent of $\varepsilon$.
By  Lemma \ref{L^pL^q} and the equality $\eqref{second order problem-B}_2$, we have
\begin{align}\label{2m base estimate-5}
&\int_{0}^T\|\Delta\tilde{c}\|^{\frac{q+2}{m}}_{L^{\frac{q+2}{m}}}\,\mathrm{d}s    \nonumber
\\
& \le  C\int_{0}^T\int_{\Omega}|(u-2\nabla g_1)\cdot\nabla \tilde{c}|^{\frac{q+2}{m}}\,\mathrm{d}x\mathrm{d}s+C\int_{0}^T\int_{\Omega}( |n\tilde{c}|^{\frac{q+2}{m}}+| \tilde{c}|^{\frac{q+2}{m}})\,\mathrm{d}x\mathrm{d}s    \nonumber
\\
& \hspace{12pt} +  C\int_{0}^T\int_{\Omega}| (|\nabla g_1|^2+\Delta g_1-u\nabla g_1-g_{1t})\tilde{c}|^{\frac{q+2}{m}}\,\mathrm{d}x\mathrm{d}s    \nonumber
\\
& \hspace{12pt} +  C\int_{0}^T\int_{\Omega}| (\Delta g_2-u\nabla g_2-ng_2-g_{2t})e^{-g_1}|^{\frac{q+2}{m}}\,\mathrm{d}x\mathrm{d}s    \nonumber
\\
& \le   C\int_{0}^T(\|u\|^{\frac{q+2}{m}}_{L^{\frac{2(q+2)}{m}}}+1)(\|\nabla\tilde{c}\|^{\frac{q+2}{m}}_{L^{\frac{2(q+2)}{m}}}+1)\,\mathrm{d}s
+C\int_{0}^T(\|u\|^{\frac{q+2}{m}}_{L^{\frac{q+2}{m}}}+\|n\|^{\frac{q+2}{m}}_{L^{\frac{q+2}{m}}})\,\mathrm{d}s    \nonumber
\\
& \le   C\int_{0}^T\|\nabla\tilde{c}\|^{\frac{2(q+2)}{m}}_{L^{\frac{2(q+2)}{m}}}\,\mathrm{d}s
+C\int_{0}^T(\|u\|^{\frac{q+2}{m}}_{L^{\frac{q+2}{m}}}+\|u\|^{\frac{2(q+2)}{m}}_{L^{\frac{2(q+2)}{m}}}+\|n\|^{\frac{q+2}{m}}_{L^{\frac{q+2}{m}}})\,\mathrm{d}s.
\end{align}

The following proof is divided into three steps.

{\bf Step 1.}\
We assume that $q\in[2m-2, 3m-2]$ such that $\frac{2(q+2)}{m}\in[4, 6]$.

By Lemma \ref{G-N}, it is easy to see that for any $\lambda>0$,
\begin{align}\label{2m base estimate-6}
\|u\|^{\frac{q+2}{m}}_{L^{\frac{q+2}{m}}}+\|u\|^{\frac{2(q+2)}{m}}_{L^{\frac{2(q+2)}{m}}}&\le C(\|u\|^{\frac{q+2}{m}}_{H^1}+\|u\|^{\frac{2(q+2)}{m}}_{H^1}),
\\\label{2m base estimate-7}
\|\nabla\tilde{c}\|^{\frac{q+2}{m}}_{L^{\frac{2(q+2)}{m}}} \le   C\|\tilde{c}\|^{\frac{q+2}{2m}}_{L^{\infty}}\|\Delta\tilde{c}\|^{\frac{q+2}{2m}}_{L^{\frac{q+2}{m}}}
+C\|\tilde{c}\|^{\frac{q+2}{m}}_{L^{\infty}}
&\le \lambda\|\Delta\tilde{c}\|^{\frac{q+2}{m}}_{L^{\frac{q+2}{m}}}+C_{\lambda}\|\tilde{c}\|^{\frac{q+2}{m}}_{L^{\infty}}.
\end{align}
Taking sufficiently small $\lambda$ and recalling Lemma \ref{L1 estimate for second order system}, we reduces \eqref{2m base estimate-5} to that
\begin{align}\label{2m base estimate-8}
&\int_{0}^T\|\Delta\tilde{c}\|^{\frac{q+2}{m}}_{L^{\frac{q+2}{m}}}\,\mathrm{d}s
\le   C\int_{0}^T\|n\|^{\frac{q+2}{m}}_{L^{\frac{q+2}{m}}}\,\mathrm{d}s+C.
\end{align}
Substituting \eqref{2m base estimate-8} into \eqref{2m base estimate-4}, we obtain
\begin{align}\label{2m base estimate-9}
&\sup_t\int_{\Omega}n^{1+q}\,\mathrm{d}x+\int_{0}^T\int_{\Omega}\left(n^{m+q-2}|\nabla n|^2+\mu n^{2+q}\right)\mathrm{d}x\mathrm{d}s
\le  C\int_{0}^T\|n\|^{\frac{q+2}{m}}_{L^{\frac{q+2}{m}}}\,\mathrm{d}s+C,
\end{align}
which implies that
\begin{align*}
\int_{0}^T\|n\|^{q+2}_{L^{q+2}}\,\mathrm{d}s
\le C\int_{0}^T\|n\|^{\frac{q+2}{m}}_{L^{\frac{q+2}{m}}}\,\mathrm{d}s+C, \quad m>1,
\end{align*}
where $q+2\in[2m, 3m]$ and $C$ is independent of $\varepsilon$, but depends on $m$, $q$, $\mu, \chi$, $\Omega$ and $T$.
After finite iterations, we will finally derive that
\begin{align}\label{2m base estimate-10}
&\sup_{t}\int_{\Omega}n^{1+q}\,\mathrm{d}x+\int_{0}^T\int_{\Omega}n^{m+q-2}|\nabla n|^{2}\,\mathrm{d}x\mathrm{d}s
+\int_{0}^T\int_{\Omega}n^{2+q}\,\mathrm{d}x\mathrm{d}s    \nonumber
\\
&\le C\int_{0}^T\|n\|^{2}_{L^2}\,\mathrm{d}s+C \le C.
\end{align}
Taking $q=3m-2$ in \eqref{2m base estimate-10} gives
\begin{align}\label{2m base estimate-11}
&\sup_{t}\int_{\Omega}n^{3m-1}\,\mathrm{d}x+\int_{0}^T\int_{\Omega}n^{4m-4}|\nabla n|^{2}\,\mathrm{d}x\mathrm{d}s
+\int_{0}^T\int_{\Omega}n^{3m}\,\mathrm{d}x\mathrm{d}s
\le C.
\end{align}

{\bf Step 2.}\
Based on \eqref{2m base estimate-11}, we prove the time-space uniform boundedness of $u$.

Note that $-\frac{3}{2(3m-1)}>-1$. Recalling \eqref{Navier-Stokes} and according to standard smoothing properties of the Stokes semigroup, there exists $\kappa>0$ such that
\begin{align*}
\|u(t)\|_{L^{\infty}}&=\int_{-\infty}^t\|e^{-(t-s)\mathcal{A}}P(n(s)\nabla\varphi(s))\|_{L^{\infty}}\,\mathrm{d}s
\\
 & \le C\int_{-\infty}^te^{-\kappa(t-s)}(t-s)^{-\frac{3}{2(3m-1)}}\|n(s)\nabla\varphi(s)\|_{L^{3m-1}}\,\mathrm{d}s
 \\
 &  \le C\int_{-\infty}^te^{-\kappa(t-s)}(t-s)^{-\frac{3}{2(3m-1)}}\|n(s)\|_{L^{3m-1}}\|\nabla\varphi(s)\|_{L^{\infty}}\,\mathrm{d}s
 \\
 &  \le C\sup_{t}\left(\|n(t)\|_{L^{3}}\|\nabla\varphi(t)\|_{L^{\infty}}\right)\int_{0}^{\infty}e^{-\kappa s}s^{-\frac{3}{2(3m-1)}}\,\mathrm{d}s
 \\
 &  \le C,
\end{align*}
which implies that
\begin{align}\label{2m base estimate-12}
\sup_t\|u(t)\|_{L^{\infty}}\le C.
\end{align}

{\bf Step 3.}\
Finally, we assume that $q\in(3m-2, \infty)$.

Combining \eqref{2m base estimate-5} with \eqref{2m base estimate-7} and \eqref{2m base estimate-12},
and taking sufficiently small $\lambda$, we have
\begin{align*}
&\int_{0}^T\|\Delta\tilde{c}\|^{\frac{q+2}{m}}_{L^{\frac{q+2}{m}}}\,\mathrm{d}s
\le   C\int_{0}^T\|n\|^{\frac{q+2}{m}}_{L^{\frac{q+2}{m}}}\,\mathrm{d}s+C.
\end{align*}
Then we can reduce \eqref{2m base estimate-4} to that
\begin{align*}
&\sup_t\int_{\Omega}n^{1+q}\,\mathrm{d}x+\int_{0}^T\int_{\Omega}\left(n^{m+q-2}|\nabla n|^2+\mu n^{2+q}\right)\mathrm{d}x\mathrm{d}s
\le  C\int_{0}^T\|n\|^{\frac{q+2}{m}}_{L^{\frac{q+2}{m}}}\,\mathrm{d}s+C.
\end{align*}
Following the same procedure as the first step, we completes the proof.
\hfill$\Box$

\vbox to 3mm{}

Next, we prove the time-space uniform boundedness of $\nabla\tilde{c}$ and $n$.
By Lemma \ref{2m base estimate}, it is easy to see that
\begin{align}\label{iteration result}
&\sup_{t}\left(\|n(t)\|_{L^4}+\|n(t)\|_{L^{2m+\frac{2}{3}}}\right) \le C,
\end{align}
where $C$ is independent of $\varepsilon$.

\begin{lemma}\label{uniform boundedness-first case}
Assume $m>1$, $\max\{2(m-1),2\}<s\le5m-1$ and $a, g, \nabla\varphi\in L^{\infty}_T(Q)$.
Let $(n, \tilde{c}, u)$ be a periodic solution of \eqref{second order problem-B}.
Then we have
\begin{align}\label{uniform boundedness-first case-1}
&\sup_{t}(\|\nabla\tilde{c} (t)\|_{L^\infty}+\|\nabla u (t)\|_{L^\infty}+\|n(t)\|_{L^\infty})\le C,
\end{align}
where $C$ is independent of $\varepsilon$.
\end{lemma}

\noindent{\bf Proof}.\
Recalling $\eqref{second order problem-B}_2$, we see that
$$
\tilde{c}_t-\Delta \tilde{c}+\tilde{c}=F(n, \tilde{c}, u),
$$
where
\begin{align*}
F(n, \tilde{c}, u)=&-(u-2\nabla g_1)\cdot\nabla \tilde{c}+(|\nabla g_1|^2+\Delta g_1+1-n-u\nabla g_1-g_{1t})\tilde{c}
\\
& +(\Delta g_2-u\nabla g_2-ng_2-g_{2t})e^{-g_1}.
\end{align*}
Notice that the time periodic solution $\tilde{c}$ of $\eqref{second order problem-B}_2$ can be expressed as follows
$$
\tilde{c}=\int_{-\infty}^te^{-(t-s)}e^{(t-s)\Delta}F\,\mathrm{d}s,
$$
where $\{e^{t\Delta}\}_{t\ge0}$ is the Neumann heat semigroup in $\Omega$,
for more properties of Neumann heat semigroup, please refer to \cite{Winkler2012}.
By lemmas \ref{L1 estimate for second order system}
and \ref{2m base estimate}, we obtain that
\begin{align*}
\|\nabla\tilde{c}(t)\|_{L^{\infty}}&=\int_{-\infty}^te^{-(t-s)}\|\nabla e^{(t-s)\Delta}F(s)\|_{L^{\infty}}\,\mathrm{d}s
\\
 & \le C\sup_{t}\|F(t)\|_{L^{4}}\int_{-\infty}^te^{-(t-s)}(t-s)^{-\frac{7}{8}}\,\mathrm{d}s
 \\
 & \le C\sup_{t}\Big((\|u\|_{L^\infty}+1)(\|\nabla \tilde{c}\|_{L^{4}}+1)+(\|n\|_{L^{4}}+1)\|\tilde{c}\|_{L^\infty}\Big)\int_{0}^{\infty}e^{-s}s^{-\frac{7}{8}}\,\mathrm{d}s
 \\
 & \le C\sup_{t}\Big(\|n\|_{L^{4}}+\|\nabla \tilde{c}\|^{\frac{1}{2}}_{L^\infty}\|\nabla \tilde{c}\|^{\frac{1}{2}}_{L^2}\Big)+C
 \\
 & \le C+C\sup_{t}\|\nabla \tilde{c}\|^{\frac{1}{2}}_{L^\infty},
\end{align*}
which implies that
\begin{align}\label{uniform boundedness-first case-3}
\sup_t\|\nabla\tilde{c}(t)\|_{L^{\infty}}\le C.
\end{align}
Similarly, we can also obtain
\begin{align}\label{uniform boundedness-first case-3'}
\sup_t\|\nabla u(t)\|_{L^{\infty}}\le C.
\end{align}

Next, we  prove the time-space uniform boundedness of $n$.
Multiplying $\eqref{second order problem-B}_1$ by $n^{q}$ with $q>3m$, and integrating it over $\Omega\times(t_0,t)$ with $t_0<t \le t_0+T$, we have
\begin{align*}
&\frac{1}{1+q}\int_{\Omega}(n^{1+q}(x,t)-n^{1+q}(x,t_0))\,\mathrm{d}x+\int_{t_0}^t\int_{\Omega}\left(qmn^{m+q-2}|\nabla n|^2+\mu n^{2+q}+ n^{1+q}\right)\,\mathrm{d}x\mathrm{d}s    \nonumber
\\
& \le q\chi\int_{t_0}^t\int_{\Omega}n^q(e^{g_1}\nabla \tilde{c}+e^{g_1}\tilde{c}\nabla g_1+\nabla g_2)\cdot\nabla n\,\mathrm{d}x\mathrm{d}s
-\chi\int_{t_0}^t\int_{\partial\Omega}e^{g_1}\tilde{c}\frac{\partial g_1}{\partial \nu}n^{1+q}\,\mathrm{d}\Gamma\mathrm{d}s    \nonumber
\\
& \hspace{12pt} +\int_{t_0}^t\int_{\Omega}((\mu an+g)n^{q}+n^{1+q})\,\mathrm{d}x\mathrm{d}s    \nonumber
\\
& \le  \frac{q}{4}\int_{t_0}^t\int_{\Omega}n^{m+q-2}|\nabla n|^2\,\mathrm{d}x\mathrm{d}s
+Cq\int_{t_0}^t\int_{\Omega}n^{q-m+2}\,\mathrm{d}x\mathrm{d}s
+C\int_{t_0}^t\int_{\partial\Omega}n^{1+q}\,\mathrm{d}\Gamma\mathrm{d}s    \nonumber
\\
& \hspace{12pt} +\frac{\mu}{4}\int_{t_0}^t\int_{\Omega}n^{2+q}\,\mathrm{d}x\mathrm{d}s+C\int_{t_0}^t\int_{\Omega}n^{q}\,\mathrm{d}x\mathrm{d}s    \nonumber
\\
& \le  \frac{q}{2}\int_{t_0}^t\int_{\Omega}n^{m+q-2}|\nabla n|^2\,\mathrm{d}x\mathrm{d}s
+Cq\int_{t_0}^t\int_{\Omega}n^{q-m+2}\,\mathrm{d}x\mathrm{d}s    \nonumber
\\
& \hspace{12pt} +\frac{\mu}{2}\int_{t_0}^t\int_{\Omega}n^{2+q}\,\mathrm{d}x\mathrm{d}s+C\int_{t_0}^t\int_{\Omega}n^{q}\,\mathrm{d}x\mathrm{d}s,
\end{align*}
which implies that
\begin{align}\label{uniform boundedness-first case-4}
&\int_{\Omega}(n^{1+q}(x,t)-n^{1+q}(x,t_0))\,\mathrm{d}x+\frac{4mq(1+q)}{(m+q)^2}\int_{t_0}^t\int_{\Omega}|\nabla n^{\frac{m+q}{2}}|^2\,\mathrm{d}x\mathrm{d}s+\int_{t_0}^t\int_{\Omega}n^{1+q}\,\mathrm{d}x\mathrm{d}s    \nonumber
\\
& \le  Cq^2\int_{t_0}^t\int_{\Omega}n^{q-m+2}\,\mathrm{d}x\mathrm{d}s+Cq\int_{t_0}^t\int_{\Omega}n^{q}\,\mathrm{d}x\mathrm{d}s.
\end{align}
By Lemma \ref{G-N} and Young inequality, for any $\lambda\in(0,1)$, we derive that
\begin{align}\label{uniform boundedness-first case-5}
&Cq^2\int_{\Omega}n^{q-m+2}\,\mathrm{d}x=Cq^2\|n^{\frac{m+q}{2}}\|_{L^{\frac{2(q-m+2)}{m+q}}}^{\frac{2(q-m+2)}{m+q}}      \nonumber
\\
& \le  Cq^2\|n^{\frac{m+q}{2}}\|_{L^{\frac{4(1+q)}{3(m+q)}}}^{\theta_1\frac{2(q-m+2)}{m+q}}\|\nabla n^{\frac{m+q}{2}}\|_{L^2}^{(1-\theta_1)\frac{2(q-m+2)}{m+q}}
+Cq^2\|n^{\frac{m+q}{2}}\|_{L^{\frac{4(1+q)}{3(m+q)}}}^{\frac{2(q-m+2)}{m+q}}      \nonumber
\\
& \le  \lambda\|\nabla n^{\frac{m+q}{2}}\|_{L^2}^{2}+Cq^{\kappa}\|n\|_{L^{\frac{2(1+q)}{3}}}^{\frac{\theta_1(m+q)(q-m+2)}{2(m-1)+\theta_1(q+2-m)}}+Cq^2\|n\|_{L^{\frac{2(1+q)}{3}}}^{q-m+2}
\end{align}
and
\begin{align}\label{uniform boundedness-first case-6}
&Cq\int_{\Omega}n^{q}\,\mathrm{d}x=Cq\|n^{\frac{m+q}{2}}\|_{L^{\frac{2q}{m+q}}}^{\frac{2q}{m+q}}      \nonumber
\\
& \le  Cq\|n^{\frac{m+q}{2}}\|_{L^{\frac{4(1+q)}{3(m+q)}}}^{\theta_2\frac{2q}{m+q}}\|\nabla n^{\frac{m+q}{2}}\|_{L^2}^{(1-\theta_2)\frac{2q}{m+q}}
+Cq\|n^{\frac{m+q}{2}}\|_{L^{\frac{4(1+q)}{3(m+q)}}}^{\frac{2q}{m+q}}      \nonumber
\\
& \le  \lambda\|\nabla n^{\frac{m+q}{2}}\|_{L^2}^{2}+Cq^{\kappa}\|n\|_{L^{\frac{2(1+q)}{3}}}^{\frac{q\theta_2(m+q)}{m+q\theta_2}}+Cq\|n\|_{L^{\frac{2(1+q)}{3}}}^{q}
\end{align}
for some $\kappa\ge2$, where $\theta_1=\frac{4(1+q)(q+2m-1)}{(7q+9m-2)(q-m+2)}$,
$\theta_2=\frac{2(1+q)(2q+3m)}{q(7q+9m-2)}$,
and the constants $C, \kappa$ are independent of $q$ and $\varepsilon$.

Substituting \eqref{uniform boundedness-first case-5}--\eqref{uniform boundedness-first case-6} into \eqref{uniform boundedness-first case-4}
and letting $\lambda$ be sufficiently small,
we obtain
\begin{align}\label{uniform boundedness-first case-7}
&\int_{\Omega}(n^{1+q}(x,t)-n^{1+q}(x,t_0))\,\mathrm{d}x+\int_{t_0}^t\int_{\Omega}n^{1+q}\,\mathrm{d}x\mathrm{d}s    \nonumber
\\
& \le  Cq^\kappa\int_{t_0}^t\|n\|_{L^{\frac{2(1+q)}{3}}}^{\frac{\theta_1(m+q)(q-m+2)}{2(m-1)+\theta_1(q+2-m)}}\mathrm{d}s+
Cq^\kappa\int_{t_0}^t\|n\|_{L^{\frac{2(1+q)}{3}}}^{\frac{q\theta_2(m+q)}{m+q\theta_2}}\mathrm{d}s      \nonumber
\\
& \hspace{12pt} +Cq^2\int_{t_0}^t\|n\|_{L^{\frac{2(1+q)}{3}}}^{q-m+2}\mathrm{d}s+Cq\int_{t_0}^t\|n\|_{L^{\frac{2(1+q)}{3}}}^{q}\mathrm{d}s.
\end{align}
Since $q\ge3m>3$, it is easy to check that
$$
\frac{2(1+q)}{3}\le\frac{\theta_1(m+q)(q-m+2)}{2(m-1)+\theta_1(q+2-m)}, \frac{q\theta_2(m+q)}{m+q\theta_2},
q-m+2, q \le 1+q
$$
Then applying Lemma \ref{periodic gronwall inequality-1} to \eqref{uniform boundedness-first case-7}, we have
\begin{align}\label{uniform boundedness-first case-8}
\sup_t\|n\|^{1+q}_{L^{1+q}}\le Cq^\kappa(\sup_t\|n\|_{L^{\frac{2(1+q)}{3}}}^{\frac{2(1+q)}{3}}+\sup_t\|n\|_{L^{\frac{2(1+q)}{3}}}^{1+q}),
\end{align}
where $C$ and $\kappa$ are independent of $q$ and $\varepsilon$.
Define a monotonically increasing sequence
\begin{align}\label{increasing sequence-1}
\{r_j\}^{\infty}_{j=0},\quad  r_{j+1}=\frac{3}{2}r_{j},\ \ r_0=2m+\frac{2}{3}.
\end{align}
Take $q=r_{j+1}-1$ in \eqref{uniform boundedness-first case-8} and $M_j=\max\{1, \sup_t\|n\|_{L^{r_j}}\}$ with $j=0,1,2,\cdots$.
Then we have
\begin{align*}
M_{j}& \le C^{\frac{1}{r_j}}(r_j-1)^{\frac{\kappa}{r_j}}(M_{j-1}^{r_{j-1}}+M_{j-1}^{r_{j}})^{\frac{1}{r_j}}
\\
&\le (2C)^{\frac{1}{r_j}}(r_j-1)^{\frac{\kappa}{r_j}} M_{j-1}
\\
& \le (2C)^{\sum_{i=1}^{j}\frac{1}{r_i}}\Pi^j_{i=1}(r_i-1)^{\frac{\kappa}{r_i}} M_{0}, \quad j=1,2,\cdots.
\end{align*}
Notice that $\sum_{i=1}^{j}\frac{1}{r_i}$ and $\Pi^{j}_{i=1}(r_i-1)^{\frac{\kappa}{r_i}}$ converge as $j\to\infty$.
Letting $j\to\infty$ gives
\begin{align}\label{uniform boundedness-first case-9}
\sup_t\|n(t)\|_{L^{\infty}}\le C+C\sup_t\|n(t)\|_{L^{2m+\frac{2}{3}}}\le C,
\end{align}
where $C$ is independent of $\varepsilon$.
The proof is completed.
\hfill$\Box$

\vbox to 3mm{}

By virtue of \eqref{2m base estimate-12} and Lemma \ref{uniform boundedness-first case}, we give
the following estimate.

\begin{lemma}\label{estimate of n-t}
Assume $m>1$, $\max\{2(m-1),2\}<s\le5m-1$, $g\ge0$  and $a, \nabla\varphi\in L^{\infty}_T(Q)$. Let $(n, \tilde{c}, u)$
be a periodic solution of \eqref{second order problem-B}.
Then we have
\begin{align}\label{estimate of n-t-1}
&\int_{0}^T\left(\|u_t\|^p_{L^p}+\|u\|^p_{W^{2,p}}+\|\tilde{c}_t\|^p_{L^p}+\|\tilde{c}\|^p_{W^{2,p}}\right)\,\mathrm{d}s\le C \ \ {\rm for\ any\ } p>1,
\\\label{estimate of n-t-2}
&\sup_t\int_{\Omega}|\nabla (n+\varepsilon)^m|^2\,\mathrm{d}x
+\int_{0}^T\int_{\Omega}(n+\varepsilon)^{m-1}\left|\frac{\partial n}{\partial t}\right|^2\,\mathrm{d}x\mathrm{d}s\le C,
\end{align}
where $C$ is independent of $\varepsilon$.
\end{lemma}

\noindent{\bf Proof}.\
Applying Lemma \ref{L^pL^q}--\ref{G-N} and the uniform boundedness of $u$, $\nabla\tilde{c}$ and $n$ to \eqref{second order problem-B}, we can obtain \eqref{estimate of n-t-1}.
Multiplying $\eqref{second order problem-B}_1$ by $\frac{1}{m}\frac{\partial(n+\varepsilon)^m}{\partial t}$, and integrating it over $\Omega\times(t_0,t)$ with $t_0<t \le t_0+T$, we have
\begin{align*}
&\frac{m}{2}\int_{\Omega}\Big((n+\varepsilon)^{2(m-1)}|\nabla n|^2(x,t)-(n+\varepsilon)^{2(m-1)}|\nabla n|^2(x,t_0)\Big)\,\mathrm{d}x
\\
& \hspace{12pt} +\int_{t_0}^t\int_{\Omega}(n+\varepsilon)^{m-1}\left|\frac{\partial n}{\partial t}\right|^2\,\mathrm{d}x\mathrm{d}s    \nonumber
\\
&  =  \int_{t_0}^t\int_{\Omega}(-u\cdot\nabla n+\mu an-\mu n^2+g)(n+\varepsilon)^{m-1}\frac{\partial n}{\partial t}\,\mathrm{d}x\mathrm{d}s    \nonumber
\\
& \hspace{12pt} -\chi\int_{t_0}^t\int_{\Omega}\nabla\cdot(e^{g_1}n\nabla \tilde{c}+e^{g_1}n\tilde{c}\nabla g_1+n\nabla g_2)(n+\varepsilon)^{m-1}\frac{\partial n}{\partial t}\,\mathrm{d}x\mathrm{d}s    \nonumber
\\
&  \le   \frac{1}{2}\int_{t_0}^t\int_{\Omega}(n+\varepsilon)^{m-1}\left|\frac{\partial n}{\partial t}\right|^2\,\mathrm{d}x\mathrm{d}s
+C\int_{t_0}^t\int_{\Omega}(n+\varepsilon)^{m-1}\left|\nabla n\right|^2\,\mathrm{d}x\mathrm{d}s    \nonumber
\\
& \hspace{12pt} +C\int_{t_0}^t\int_{\Omega}\left|\Delta\tilde{c}\right|^2\,\mathrm{d}x\mathrm{d}s+C,
\end{align*}
which implies that
\begin{align*}
&\int_{\Omega}\Big((n+\varepsilon)^{2(m-1)}|\nabla n|^2(x,t)-(n+\varepsilon)^{2(m-1)}|\nabla n|^2(x,t_0)\Big)\,\mathrm{d}x    \nonumber
\\
& \hspace{12pt} +\int_{t_0}^t\int_{\Omega}(n+\varepsilon)^{m-1}\left|\frac{\partial n}{\partial t}\right|^2\,\mathrm{d}x\mathrm{d}s    \nonumber
\\
&  \le   C\int_{t_0}^t\int_{\Omega}(n+\varepsilon)^{m-1}\left|\nabla n\right|^2\,\mathrm{d}x\mathrm{d}s
+C\int_{t_0}^t\int_{\Omega}\left|\Delta\tilde{c}\right|^2\,\mathrm{d}x\mathrm{d}s+C.
\end{align*}
Combining with Lemma \ref{periodic gronwall inequality-1}, \ref{L1 estimate for second order system}, and
\ref{uniform boundedness-first case},
we derive that
\begin{align*}
&\sup_t\int_{\Omega}(n+\varepsilon)^{2(m-1)}|\nabla n|^2(x,t)\,\mathrm{d}x
+\int_{0}^T\int_{\Omega}(n+\varepsilon)^{m-1}\left|\frac{\partial n}{\partial t}\right|^2\,\mathrm{d}x\mathrm{d}s    \nonumber
\\
&  \le   C\int_{0}^T\int_{\Omega}(n+\varepsilon)^{m-1}\left|\nabla n\right|^2\,\mathrm{d}x\mathrm{d}s
+C\int_{0}^T\int_{\Omega}\left|\Delta\tilde{c}\right|^2\,\mathrm{d}x\mathrm{d}s+C    \nonumber
\\
&  \le   C\int_{0}^T\int_{\Omega}(n+\varepsilon)^{m-1}|\nabla \sqrt{n}|^2\,\mathrm{d}x\mathrm{d}s+C
\\
& \le C,
\end{align*}
where $C$ is independent of $\varepsilon$. The proof is completed.
\hfill$\Box$

\vbox to 3mm{}

The proof of Proposition \ref{estimate for second order system} is a consequence of Lemma \ref{L1 estimate for second order system}--\ref{estimate of n-t}.

\noindent{\bf Proof of Theorem \ref{main result-1}}.\
Let $(u_{\varepsilon}, \tilde{c}_{\varepsilon}, n_{\varepsilon})$
be a time periodic solution of the problem \eqref{second order problem-B} satisfying Proposition \ref{estimate for second order system}.
Then we have
\begin{align}
&-\iint_{Q_T}n_{\varepsilon}\phi_{1t}\,\mathrm{d}x\mathrm{d}s+\iint_{Q_T}\nabla(n_{\varepsilon}+\varepsilon)^m\nabla\phi_1\,\mathrm{d}x\mathrm{d}s
+\iint_{Q_T}(\varepsilon n_{\varepsilon}^{s+1}+u_{\varepsilon}\cdot\nabla n_{\varepsilon})\phi_1\,\mathrm{d}x\mathrm{d}s      \nonumber
\\
&  \hspace{12pt} =   \chi\iint_{Q_T}(e^{g_1}n_{\varepsilon}\nabla \tilde{c}_{\varepsilon}+n_{\varepsilon}\tilde{c}_{\varepsilon}\nabla e^{g_1}+n_{\varepsilon}\nabla g_2)\cdot\nabla\phi_1\,\mathrm{d}x\mathrm{d}s
-\chi\int_0^T\int_{\partial\Omega}n_{\varepsilon}\tilde{c}_{\varepsilon}\phi_1\frac{\partial e^{g_1}}{\partial \nu}\,\mathrm{d}\Gamma\mathrm{d}s      \nonumber
\\\label{approximation of weak solution-n}
& \hspace{24pt} +\iint_{Q_T}(\mu n_{\varepsilon}(a-n_{\varepsilon})+g)\phi_1\,\mathrm{d}x\mathrm{d}s,
\\[2mm]
&-\iint_{Q_T}\tilde{c}_{\varepsilon}\phi_{2t}\,\mathrm{d}x\mathrm{d}s+\iint_{Q_T}\nabla\tilde{c}_{\varepsilon}\nabla\phi_2\,\mathrm{d}x\mathrm{d}s
+\iint_{Q_T}((u_{\varepsilon}-2\nabla g_1)\cdot\nabla \tilde{c}_{\varepsilon})\phi_2\,\mathrm{d}x\mathrm{d}s      \nonumber
\\
&\hspace{12pt} =   \iint_{Q_T}(|\nabla g_1|^2+\Delta g_1-n_{\varepsilon}-u_{\varepsilon}\nabla g_1-g_{1t})\tilde{c}_{\varepsilon}\phi_2\,\mathrm{d}x\mathrm{d}s      \nonumber
\\\label{approximation of weak solution-c}
&\hspace{24pt} +\iint_{Q_T}(\Delta g_2-u_{\varepsilon}\nabla g_2-n_{\varepsilon}g_2-g_{2t})e^{-g_1}\phi_2\,\mathrm{d}x\mathrm{d}s,
\\[2mm]\label{approximation of weak solution-u}
&-\iint_{Q_T}u_{\varepsilon}\phi_{3t}\,\mathrm{d}x\mathrm{d}s+\iint_{Q_T}\nabla u_{\varepsilon}\nabla\phi_3\,\mathrm{d}x\mathrm{d}s
=\iint_{Q_T}n_{\varepsilon}\nabla \varphi\phi_3\,\mathrm{d}x\mathrm{d}s,
\end{align}
for any $\phi_{1}, \phi_{2}, \phi_{3}\in H^1_T(Q)$ with $\frac{\partial \phi_{1,2}}{\partial \nu}|_{\partial\Omega}=0$, $\phi_3|_{\partial\Omega}=0$ and $\nabla\cdot\phi_3=0$.
Using Sobolev imbedding theorem and taking $\varepsilon\to0$, we have (if necessary, we may choose a subsequence)
\begin{equation*}
\begin{array}{l@{\hspace{-1pt}}l@{\hspace{30pt}}l}
&u_{\varepsilon}\to u,  \ \tilde{c}_{\varepsilon}\to \tilde{c},  &  {\rm uniformly},
\\[2mm]
&u_{\varepsilon}\rightharpoonup u,  \ \tilde{c}_{\varepsilon}\rightharpoonup \tilde{c}, &  {\rm in}\ \  W_p^{2,1}(Q_T)\ \ {\rm for\ any}\ \ p>1,
\\[2mm]
&n_{\varepsilon},\ n_{\varepsilon}+\varepsilon\to n, & {\rm in}\ \ L^p(Q_T)\ \ {\rm for\ any}\ \ p>1,
\\[2mm]
&\varepsilon n^{s+1}_{\varepsilon}\to 0,  \ n_{\varepsilon}\stackrel{^*}{\rightharpoonup} n, & {\rm in}\ \ L^\infty(Q_T),
\\[2mm]
&\nabla(n_{\varepsilon}+\varepsilon)^m\rightharpoonup\nabla n^m,  & {\rm in}\ \  L^2(Q_T).
 \end{array}
\end{equation*}
Then $(u, \tilde{c}, n)$ is a time periodic solution of the problem \eqref{homogeneous problem} satisfies
\eqref{main result-1-1}--\eqref{main result-1-3}.
\hfill$\Box$

\vspace{5mm}

\setcounter{equation}{0}
\section{Strong time periodic solutions}

In this section, we improve the regularity for $m\in(1,\frac{4}{3}]$ and prove that the obtained
time periodic solution is strong solution. Let $(n_{\varepsilon}, \tilde{c}_{\varepsilon}, u_{\varepsilon})$
be a time periodic solution of \eqref{second order problem-B}.
We still assume that $(n_{\varepsilon}, \tilde{c}_{\varepsilon}, u_{\varepsilon})$ is sufficiently
smooth. Otherwise, the following estimates can be obtained by an approximating
process.

\begin{lemma}\label{estimate1 of strong solution}
Assume $1<m\le\frac{4}{3}$, $\max\{2(m-1),2\}<s\le5m-1$, $g\ge0$  and $a, \nabla\varphi\in L^{\infty}_T(Q)$. Let $(n_{\varepsilon}, \tilde{c}_{\varepsilon}, u_{\varepsilon})$
be a time periodic solution of \eqref{second order problem-B}.
Then we have
\begin{align}\label{estimate1 of strong solution-1}
&\sup_t\int_{\Omega}|\nabla \sqrt{n_{\varepsilon}+\varepsilon}|^2\,\mathrm{d}x
+\int_{0}^T\int_{\Omega}(n_{\varepsilon}+\varepsilon)^{m-4}|\nabla n_{\varepsilon}|^{4}\,\mathrm{d}x\mathrm{d}s\le C,
\end{align}
where $C$ is independent of $\varepsilon$.
\end{lemma}

\noindent{\bf Proof}.\
The proof is the similar with that of Lemma 3.5 in \cite{Jin2021arXiv}. Actually, the proof is more concise in our case
by applying the estimate \eqref{estimate of n-t-1} and the uniform boundedness of $u_{\varepsilon}$, $\tilde{c}_{\varepsilon}$ and $n_{\varepsilon}$.
So we omit it here.
\hfill$\Box$

\vbox to 3mm{}

\begin{lemma}\label{estimate2 of strong solution}
Assume $1<m\le\frac{4}{3}$, $\max\{2(m-1),2\}<s\le5m-1$, $g\ge0$  and $a, \nabla\varphi\in L^{\infty}_T(Q)$. Let $(n_{\varepsilon}, \tilde{c}_{\varepsilon}, u_{\varepsilon})$
be a time periodic solution of \eqref{second order problem-B}.
Then we have
\begin{align}\label{estimate2 of strong solution-1}
&\int_{0}^T\Big(\|n_{\varepsilon t}\|_{L^2}^2+\|\Delta(n_{\varepsilon}+\varepsilon)^m\|_{L^2}^2\Big)\,\mathrm{d}s\le C,
\end{align}
where $C$ is independent of $\varepsilon$.
\end{lemma}

\noindent{\bf Proof}.\
Multiplying $\eqref{second order problem-B}_1$ by $n_{\varepsilon t}$, integrating it over $\Omega\times(t_0,t)$ with $t_0<t \le t_0+T$,
and applying \eqref{estimate of n-t-1}, \eqref{estimate1 of strong solution-1} and the uniform boundedness of $u_{\varepsilon}$, $\tilde{c}_{\varepsilon }$ and $n_{\varepsilon }$, we have
\begin{align*}
&\frac{m}{2}\int_{\Omega}\Big((n_{\varepsilon}+\varepsilon)^{m-1}|\nabla n_{\varepsilon}|^2(x,t)-(n_{\varepsilon}+\varepsilon)^{m-1}|\nabla n_{\varepsilon}|^2(x,t_0)\Big)\,\mathrm{d}x+\int_{t_0}^t\int_{\Omega}\left|n_{\varepsilon t}\right|^2\,\mathrm{d}x\mathrm{d}s
\\
&  =  -\chi\int_{t_0}^t\int_{\Omega}\nabla\cdot(e^{g_1}n_{\varepsilon}\nabla \tilde{c}_{\varepsilon}+e^{g_1}n_{\varepsilon}\tilde{c}_{\varepsilon}\nabla g_1+n_{\varepsilon}\nabla g_2)n_{\varepsilon t}\,\mathrm{d}x\mathrm{d}s    \nonumber
\\
& \hspace{12pt}   +\int_{t_0}^t\int_{\Omega}(-\varepsilon n_{\varepsilon}^{s+1}-u_{\varepsilon}\cdot\nabla n_{\varepsilon}+\mu an_{\varepsilon}-\mu n_{\varepsilon}^2+g)n_{\varepsilon t}\,\mathrm{d}x\mathrm{d}s  \nonumber
\\
& \hspace{12pt}   +\frac{m(m-1)}{2}\int_{t_0}^t\int_{\Omega}(n_{\varepsilon}+\varepsilon)^{m-2}n_{\varepsilon t}|\nabla n_{\varepsilon}|^2\,\mathrm{d}x\mathrm{d}s  \nonumber
\\
&  \le   \frac{1}{2}\int_{t_0}^t\int_{\Omega}\left|n_{\varepsilon t}\right|^2\,\mathrm{d}x\mathrm{d}s
+C\int_{t_0}^t\int_{\Omega}(n_{\varepsilon}+\varepsilon)^{2m-4}\left|\nabla n_{\varepsilon}\right|^4\,\mathrm{d}x\mathrm{d}s+C    \nonumber
\\
&  \le   \frac{1}{2}\int_{t_0}^t\int_{\Omega}\left|n_{\varepsilon t}\right|^2\,\mathrm{d}x\mathrm{d}s
+C\int_{t_0}^t\int_{\Omega}(n_{\varepsilon}+\varepsilon)^{m-4}\left|\nabla n_{\varepsilon}\right|^4\,\mathrm{d}x\mathrm{d}s+C    \nonumber
\\
&  \le   \frac{1}{2}\int_{t_0}^t\int_{\Omega}\left|n_{\varepsilon t}\right|^2\,\mathrm{d}x\mathrm{d}s+C.
\end{align*}
Combining with Lemma \ref{periodic gronwall inequality-1},
we derive that
\begin{align}\label{estimate2 of strong solution-2}
&\int_{0}^T\int_{\Omega}\left|n_{\varepsilon t}\right|^2\,\mathrm{d}x\mathrm{d}s\le C,
\end{align}
where $C$ is independent of $\varepsilon$.
Similarly, multiplying $\eqref{second order problem-B}_1$ by $-\Delta(n_{\varepsilon}+\varepsilon)^m$, integrating it over $\Omega\times(t_0,t)$ with $t_0<t \le t_0+T$,
and applying \eqref{estimate of n-t-1}, \eqref{estimate1 of strong solution-1}, \eqref{estimate2 of strong solution-2} and the uniform boundedness of $u_{\varepsilon}$, $\tilde{c}_{\varepsilon}$ and $n_{\varepsilon}$,
we also obtain that
\begin{align}\label{estimate2 of strong solution-3}
&\int_{0}^T\int_{\Omega}\left|\Delta(n_{\varepsilon}+\varepsilon)^m\right|^2\,\mathrm{d}x\mathrm{d}s\le C,
\end{align}
where $C$ is independent of $\varepsilon$.
The proof is completed.
\hfill$\Box$

\vbox to 3mm{}

Next, we prove Theorem \ref{main result-2}.

\noindent{\bf Proof of Theorem \ref{main result-2}}.\
Let $(u_{\varepsilon}, \tilde{c}_{\varepsilon}, n_{\varepsilon})$
be a weak time periodic solution of the problem \eqref{second order problem-B} satisfying Proposition \ref{estimate for second order system}.
Using Lemma \ref{estimate2 of strong solution} and Sobolev imbedding theorem
and taking $\varepsilon\to0$, we obtain that (if necessary, we may choose a subsequence)
\begin{equation*}
\begin{array}{l@{\hspace{-1pt}}l@{\hspace{30pt}}l}
&u_{\varepsilon}\to u,  \ \tilde{c}_{\varepsilon}\to \tilde{c},  &  {\rm uniformly},
\\[2mm]
&u_{\varepsilon}\to u,  \ \tilde{c}_{\varepsilon}\to \tilde{c}, & {\rm in}\ \  L^p(Q_T)\ \ {\rm for\ any}\ \ p>1,
\\[2mm]
&\tilde{c}_{\varepsilon t}\rightharpoonup \tilde{c}_{t},  \ \nabla\tilde{c}_{\varepsilon}\to \nabla\tilde{c},\ \ \Delta\tilde{c}_{\varepsilon}\rightharpoonup\Delta\tilde{c},  & {\rm in}\ \  L^p(Q_T)\ \ {\rm for\ any}\ \ p>1,
\\[2mm]
&u_{\varepsilon t}\rightharpoonup u_{t},  \ \nabla u_{\varepsilon}\to\nabla u, \ \ \Delta u_{\varepsilon}\rightharpoonup\Delta u,  & {\rm in}\ \  L^p(Q_T)\ \ {\rm for\ any}\ \ p>1,
\\[2mm]
&n_{\varepsilon}, n_{\varepsilon}+\varepsilon\to n, & {\rm in}\ \ L^p(Q_T)\ \ {\rm for\ any}\ \ p>1,
\\[2mm]
&\varepsilon n^{s+1}_{\varepsilon}\to 0,  \ n_{\varepsilon}\stackrel{^*}{\rightharpoonup} n, & {\rm in}\ \ L^\infty(Q_T),
\\[2mm]
&\nabla n_{\varepsilon}\to \nabla n, & {\rm in}\ \ L^p(Q_T)\ \ {\rm for\ any}\ \ p\in(1,6),
\\[2mm]
&n_{\varepsilon t}\rightharpoonup n_{t},  \ \Delta(n_{\varepsilon}+\varepsilon)^m\rightharpoonup\Delta n^m,  & {\rm in}\ \  L^2(Q_T).
 \end{array}
\end{equation*}
Applying integration by parts to \eqref{approximation of weak solution-n}--\eqref{approximation of weak solution-u} and letting $\varepsilon\to0$,
we derive that
\begin{align}
&\iint_{Q_T}(n_t-\Delta n^m+u\cdot\nabla n)\phi_1\,\mathrm{d}x\mathrm{d}s-\iint_{Q_T}(\mu n_{\varepsilon}(a-n_{\varepsilon})+g)\phi_1\,\mathrm{d}x\mathrm{d}s      \nonumber
\\\label{approximation of strong solution-n}
&\hspace{12pt} =   -\chi\iint_{Q_T}\nabla\cdot(e^{g_1}n_{\varepsilon}\nabla \tilde{c}_{\varepsilon}+n_{\varepsilon}\tilde{c}_{\varepsilon}\nabla e^{g_1}+n_{\varepsilon}\nabla g_2)\phi_1\,\mathrm{d}x\mathrm{d}s,
\\[2mm]
&\iint_{Q_T}(\tilde{c}_t-\Delta \tilde{c}+(u-2\nabla g_1)\cdot\nabla \tilde{c})\phi_2\,\mathrm{d}x\mathrm{d}s      \nonumber
\\
& \hspace{12pt} =   \iint_{Q_T}(|\nabla g_1|^2+\Delta g_1-n-u\nabla g_1-g_{1t})\tilde{c}\phi_2\,\mathrm{d}x\mathrm{d}s        \nonumber
\\\label{approximation of strong solution-c}
&\hspace{24pt} +\iint_{Q_T}(\Delta g_2-u\nabla g_2-ng_2-g_{2t})e^{-g_1}\phi_2\,\mathrm{d}x\mathrm{d}s,
\\[2mm]\label{approximation of strong solution-u}
&\iint_{Q_T}(u_t-\Delta u  +\nabla \pi -n\nabla \varphi)\phi_3\,\mathrm{d}x\mathrm{d}s
=0,
\end{align}
for any $\phi_{1}, \phi_{2}, \phi_{3}\in H^1_T(Q)$ with $\frac{\partial \phi_{1,2}}{\partial \nu}|_{\partial\Omega}=0$, $\phi_3|_{\partial\Omega}=0$ and $\nabla\cdot\phi_3=0$.
Then $(u, \tilde{c}, n)$ is a strong time periodic solution of the problem \eqref{homogeneous problem} satisfying Proposition \ref{estimate for second order system}
and Lemma \ref{estimate2 of strong solution}. The proof is completed.
\hfill$\Box$

\vspace{5mm}

\section*{Acknowledgements}

The authors would like to thank the anonymous referee for the valuable comments
and suggestions.

\vspace{5mm}

%-------------------------------------------------------------------

\end{document}